

\magnification 1200
\input amssym.def
\input amssym
\catcode`\@=12
\font\smalltype=cmr10 at 10truept
 at 10truept
\nopagenumbers

\def\leftheadline{\ifnum\pageno >1
	\tenrm \folio\hfil \smalltype
	 Bradlow, Daskalopoulos, and Wentworth
	\hfil\fi}

\def\rightheadline{\ifnum\pageno >1
	\hfil \smalltype
	Birational Equivalences of Vortex Moduli
	\hfil \tenrm \folio\fi}

\headline={\ifodd\pageno\rightheadline\else\leftheadline\fi}

$\hbox{}$
\bigskip\bigskip
\centerline{\bf BIRATIONAL EQUIVALENCES OF VORTEX MODULI
	\footnote{${}^\star$}{\smalltype Revised -- April 6, 1993.}
	}
\vskip .75in
\centerline{Steven B. Bradlow \footnote{${}^1$}{\smalltype
Supported in part by NSF grant DMS-9103950 and an NSF-NATO
Postdoctoral Fellowship.}
}
\medskip
\centerline{Georgios D. Daskalopoulos}
\medskip
\centerline{Richard A. Wentworth \footnote{${}^2$}{\smalltype
Supported in part by NSF Mathematics Postdoctoral Fellowship
DMS-9007255.}}
\bigskip\bigskip
\midinsert
\narrower\narrower \noindent
{\bf Abstract.}$\quad$
We construct a finite dimensional
K\"ahler manifold with a holomorphic, symplectic
circle action whose symplectic reduced spaces may be identified with
the $\tau$-vortex moduli spaces (or $\tau$-stable pairs).   The
Morse theory of the circle action induces  natural birational maps
between the reduced spaces for different values of $\tau$ which in
the case of rank two bundles can be canonically resolved in a
sequence of blow-ups and blow-downs.
\endinsert
\eject

\def \aspace{{\cal C} \times \kform 0E}

\def \Bhat{\hat{\cal B}}
\def \Btau{{\cal B}_{\tau}}
\def \curvdh{F_{\dbare,H}}
\def \Cx{C^{\dbare}_{\phi}}
\def \Cxo{C^{\dbare}_{\phi,0}}

\def \dbare{\bar{\partial}_E}
\def \dbf{\dbare,\phi}

\def \G{\frak G}

\def \Gc{{\frak G}^{\Bbb C}}
\def \Gci{{\frak G}_1^{\Bbb C}}
\def \Gco{{\frak G}^{\Bbb C}_0}
\def \H{{\cal H}}

\def \kform#1#2{\Omega^{#1}(#2)}

\def \mapright#1{\smash{\mathop{\longrightarrow}
	\limits^{#1}}}

\def \mmap{\Psi}
\def \mmapo{\Psi_0}
\def \muM{\mu_M}
\def \mum{\mu_m (\phi)}
\def \pqform#1#2#3{\Omega^{#1,#2}(#3)}

\def \ss{{\cal C}_{ss}}

\def \tauhat{{\tau\hbox{\smalltype Vol}(\Sigma)\over 4\pi}}
\def \Tr{\hbox{Tr}\,}
\def \U{\hbox{U}(1)}

\def \Vt{{\cal V}_{\tau}}

\def\BtL{\Btau(L)}
\def\Jd{{{\cal J}_d}}

\def\CBbb{{\Bbb C}}
\def\ZBbb{{\Bbb Z}}
\def\degree{\hbox{degree}}
\def\End{\hbox{End}\,}

\def\Lie{\hbox{Lie}\,}
\def\Rank{\hbox{Rank}}
\def\rank{\hbox{rank}}

\def\Wtau{{\cal W}_\tau}
\def\Ztau{{\cal Z}_\tau}
\def\tuple{(d_\phi, R_\phi, R_1,\ldots, R_n)}
\def\lra{\longrightarrow}
\def\ad{\hbox{ad}\,}

\baselineskip 16pt \normalbaselineskip 16pt

\beginsection 1. Introduction

For holomorphic bundles over a Riemann surface there is essentially
one notion of stability, and hence a single moduli space for bundles
of fixed rank and degree.
This rigidity can disappear when one considers moduli of bundles over
higher dimensional varieties or when one considers bundles with
additional structure, such as parabolic bundles.  The concept of
stability can then depend on parameters, and one can get families of
moduli .  In this paper we explore this phenomenon in the case of
holomorphic bundles with prescribed global sections --- the so-called
holomorphic pairs.  The point of view we take is inspired by Morse
Theory and symplectic geometry.

In [B1] and  [B-D1]
we introduced a notion of stability for a pair
$(E,\phi)$\ consisting of a holomorphic bundle together with a
holomorphic section. The definition involves a real valued parameter
and can be stated as follows:

\proclaim Definition 1.1.  Let $E\longrightarrow \Sigma$ be a
holomorphic vector bundle over a compact
Riemann surface $\Sigma$. Let $\phi\in
H^0(\Sigma ,E)$ be a holomorphic section,  and
let $\tau$\ be a real number.
We say that the pair $(E,\phi)$  is $\tau$-stable (resp.
 $\tau$-semistable) if the following
two conditions hold:
\itemitem{(i)} ${\hbox{\smalltype degree}(F)\over
\hbox{\smalltype rank}(F)}<\tau$ (resp. $\leq \tau$),
for every holomorphic
subbundle $F\subset E$;
\itemitem{(ii)} ${\hbox{\smalltype degree}(E/F)\over
\hbox{\smalltype rank}(E/F)}>\tau$ (resp. $\geq \tau$),
for every proper
holomorphic subbundle $F\subset E$ such that $\phi$ is a section of
$F$.

{\it Throughout the paper}, we shall denote the rank of $E$ by $R$,
the
degree of $E$ by $d$, and the genus of $\Sigma$ by $g$.  We shall
also assume that $g\geq 2$, $R\geq 2$,
and that $d > R(2g-2)$ (cf. Assumption 2
of [B-D1]).
These assumptions may be relaxed, giving rise to
interesting special cases; for the presentation of the general
theory, however, it is convenient to make them.

Definition 1.1, and specifically the origin of the parameter $\tau$,
is motivated by a correspondence between stability criteria and the
existence of special bundle metrics.  In the case of pure holomorphic
bundles, this is the Hitchin-Kobayashi correspondence between
stability and the Hermitian-Einstein condition [Ko].  The
Hermitian-Einstein condition is expressed in the form of a set of
partial differential equations called the Hermitian-Einstein
equations.  These amount to a pointwise constraint on the curvature
of a unitary connection.  But curvature forms are always globally
constrained, via Chern-Weil formulae, by the topology of the bundle.
For bundles over closed Riemann surfaces, this removes any ambiguity
in the Hermitian-Einstein condition and hence in the definition of
stability. Now the Hermitian-Einstein equations admit a natural
modification which is appropriate when a global section is
prescribed [B1]. The new equations, called the vortex equations, are
obtained by adding extra terms which involve only the global section.
These terms are not subject to any topological constraint, and thus
unlike the Hermitian-Einstein equations, the vortex equations involve
a true parameter. Tracing back the Hitchin-Kobayashi correspondence
from the equations to a constraint on the holomorphic structure, one
is led from the vortex equations to the above notion of stability.
Since the equations have a parameter, so does the notion of
stability.

The parameter $\tau$
can be explained in another way, which depends on a
correspondence between holomorphic pairs on a compact Riemann surface
$\Sigma$ and a certain holomorphic extension on $\Sigma\times\Bbb
P^1$. This correspondence was observed by Garcia-Prada, who used it
to relate the stability of a pair to the stability of the
corresponding extension [G-P]. But on $\Sigma\times\Bbb P^1$, the
definition of stability depends on a choice of polarization. This
introduces a single parameter which is essentially the relative
weights of the polarizations on $\Sigma$ and $\Bbb P^1$. When
transferred back to the pair on $\Sigma$, the parameter is no longer
in the polarization, but in the stability criterion itself.

The impact of the parameter $\tau$ is shaped primarily by two
things.  Firstly, for purely numerical reasons, at almost all values
of $\tau$, the strict inequalities in Definition 1.1 are equivalent to
weak inequalities.  At only a discrete set of values (specifically,
rational numbers whose denominator is strictly
between 0 and $\Rank(E)$) is
equality possible.  Let us call these values the critical values of
$\tau$.  Secondly, for values of $\tau$ between any two successive
critical values, the definitions of stability are entirely
equivalent.  Furthermore, for these ``generic" values of $\tau$
we get good moduli spaces of $\tau$-stable pairs. Specifically, we
have

\proclaim  Theorem 1.2. {\tenrm (cf. [B-D1], [B-D2], [Be], [G-P], [Th])}
Let $\Btau$ denote the set of isomorphism classes of $\tau$-stable
pairs on $E$. If $\tau$ is not a
positive rational with denominator less than $R$, then $\Btau$
naturally has the structure of a compact K\"ahler manifold of
dimension $d+(R^2-R)(g-1)$.  Indeed, $\Btau$ is an algebraic
variety (see Theorem 6.3),
and the same result holds for $\Btau(L)$, where
$L\to\Sigma$ is a degree $d$ line bundle and $\Btau(L)$ denotes the
space of $\tau$-stable pairs with fixed determinant $L$.

If we think of the parameter as a ``height function" and the
moduli spaces as ``level sets", then these features strongly
suggest a Morse Theory interpretation. After all, the range of the
height function is partitioned into intervals between critical
values, and the levels sets at levels strictly between successive
critical heights are all equivalent.

In this paper we will show that this is more than simply an analogy.
We will give a precise way of realizing just such a picture.  In fact
the parameter $\tau$ can be realized as a height with respect to a
Morse function of a very special kind, namely one arising from a
symplectic moment map.  This is not too surprising since it is well
known that the equations corresponding to stability criteria (i.e.
the Hermitian-Einstein and Vortex equations) have moment map
interpretations.  By using this aspect of the problem, we relate  the
present situation to a phenomenon studied in symplectic geometry,
i.e. the variation of symplectically reduced level sets of moment
maps [G-S].

In essence what we do is to construct a single large ``master
space" (the terminology is due to Bertram)
which contains the stable and semistable pairs for all values
of the parameter $\tau$.  The space has a symplectic structure and a
symplectic  circle action.  We detect $\tau$ as the value of the
moment map for this circle action, and we recover the moduli spaces
of $\tau$-stable pairs as the Marsden-Weinstein reductions for
different values of this moment map.
Stated more precisely, we prove the following

\proclaim Theorem 1.3.  Consider the holomorphic pairs on $E$.
There is a compact topological space $\Bhat$ whose points correspond
to holomorphic pairs which are $\tau$-semistable for  at least one
value of $\tau$.  Furthermore, there is an open set
$\Bhat_0\subset\Bhat$ which has a natural K\"ahler manifold
structure.  The space $\Bhat$ and $\Bhat_0$ have the following
properties:
\itemitem{(i)}  There is a quasi-free $\U$-action on $\Bhat$, i.e.
an action for which the isotropy subgroup is either trivial or the
whole $\U$.  On $\Bhat_0$ the action is holomorphic and symplectic.
\itemitem{(ii)}  There is a moment map $f:\Bhat_0\to{\Bbb R}$ for
this $\U$-action which extends continuously to $\Bhat$.  The
critical values for $f$ are precisely the critical values of
the parameter $\tau$.
\itemitem{(iii)} The level sets $f^{-1}(\tau)$ are $\U$-invariant.
At regular values, the orbit spaces $f^{-1}(\tau)/\U$ inherit  a
K\"ahler structure and can be identified with the moduli spaces
$\Btau$.  At critical values the orbit spaces correspond to the
spaces of isomorphism classes of semistable pairs.

\noindent{\sl In the case of rank two bundles of odd degree, $\Bhat$
itself has the structure of a K\"ahler V-manifold with
at most ${\Bbb Z}_2$ singularities
 along the minimum value of $f$. Furthermore, $f$ is a perfect Morse
function in the sense of Bott.}

An important feature of our picture is that the master space $\Bhat$
is  compact, and thus the critical values of the $\U$ moment map
include an absolute maximum and an absolute minimum.  The
level sets at these extremal values correspond to moduli spaces of
semistable {\it bundles}. If  the degree $d$ of the underlying bundle
$E$ is coprime to both $R$ and
$(R-1)$, then at one extreme we obtain precisely the moduli space of
rank $R$ degree $d$
stable bundles, while at the other extreme we obtain
precisely the moduli space of rank $(R-1)$ degree $d$ stable bundles.

In the case of rank two
bundles we can use our construction to recover
some of the beautiful results of Thaddeus [T] (see also
Theorem 6.4). Using  techniques from Geometric Invariant Theory,
Thaddeus constructed and
analyzed the moduli spaces of $\tau$-stable pairs
with fixed determinant and rank two.
He showed that  for values of $\tau$
separated by a critical value, the moduli spaces are related by flip
in the sense of Mori theory.  That is, the spaces are birationally
equivalent projective varieties, and yield the same space if each is
blown up along the locus where it differs from the other.  Thus the
one is transformed into the other by first blowing up, and then
blowing down the exceptional divisor ``along a different
direction".  In our master space construction this phenomenon has an
explanation both from the symplectic point of view as well as in
terms of the Morse theory.  From the symplectic point of view it
corresponds exactly to the relationship between reduced level sets of
moment maps as described by Guillemin and Sternberg in [G-S].  In
terms of the Morse theory, the birationality of the level sets comes
from a map induced by flows along the gradient lines of the Morse
function. The centers of the blow-up
in the two spaces are given by the
stable and unstable manifolds in the sense of Morse theory (i.e. the
points on flow lines which terminate at critical points).
Furthermore, by again using the gradient flow, the
exceptional divisors of both blow-ups can be
identified as the projectivized normal bundle of the critical
submanifold in the critical level set.

\noindent {\it Acknowledgement.}  The authors are pleased to
acknowledge the warm hospitality of the Mathematics Institute at the
University of Warwick where part of this work was completed.

\beginsection 2.  Moment maps and master spaces

\centerline{\it \S 2.1 Outline of the construction}

In this section we carry out the construction of the
space $\Bhat$ described
in Theorem
1.3. The construction is formally very similar to that given in
[B-D1] for the moduli spaces
of $\tau$-stable pairs .  We begin with a brief
overview of both the construction
of the $\Btau$\ and the modifications
required for the new space $\Bhat$.

As in the  Introduction, let $E\to\Sigma$ be a fixed complex vector
bundle of rank $R$ and degree $d$.
Also (cf. [B-D1] for more details) let
${\cal C}$ denote the space
of $\bar\partial$-operators on $E$ (or equivalently, the
space of holomorphic structures on $E$), and let $\Omega^0(E)$
denote the space of smooth sections of $E$.
The space of holomorphic pairs is
then given by
$${\cal H}=\left\{ (\dbf)\in {\cal C}\times\Omega^0(E) :
\dbare\phi=0\right\}\leqno(2.1)$$
(cf. [B-D1], Definition 1.1, and note that here we allow the case
$\phi\equiv 0$).
The complex gauge group $\Gc$, i.e. the group of bundle
automorphisms, acts on $\cal H$\ by
$$g(\dbf)=(g\circ\dbare\circ g^{-1}, g\phi).\leqno(2.2)$$
The $\Gc$-orbits correspond to
isomorphism classes of holomorphic pairs.  For  the
construction of the moduli spaces $\Btau$, we need
to identify the orbits
corresponding to the $\tau$-stable pairs. We use

\proclaim Theorem 2.1. {\tenrm (see [B1])}  Let $E\to\Sigma$ be a
fixed complex vector bundle over a closed Riemann surface, and let
$(\dbare,\phi)$ be a holomorphic pair as in Definition 1.1.
Suppose that $(\dbf)$
is $\tau$-stable for a given value of the parameter
$\tau$. Then the $\tau$-Vortex equation
$$\sqrt{-1}\Lambda\curvdh +{1\over 2}\phi\otimes\phi^\ast={\tau\over
2}{\bf I}\leqno(2.3)$$
considered as an equation for the metric $H$,
has a unique smooth solution.
Here $\curvdh$ is the curvature
of a metric connection, $\Lambda\curvdh$ is a
section in $\kform 0 {\End E}$ and is
obtained by a contraction of $\curvdh$
against the K\"ahler form on $\Sigma$,
$\phi\otimes\phi^\ast$\ is a section of
$\kform 0 {E\otimes E^*}\simeq\kform 0 {\End E}$,
and ${\bf I}$\ is the identity
section in $\kform 0 {\End E}$.
Conversely, suppose that for a given value of $\tau$ there is a
Hermitian metric $H$ on $E$ such that the $\tau$-vortex equation
is satisfied by $(\dbf, H)$.  Then $E$ splits holomorphically as
$E=E_\phi\oplus E_s$, where
\itemitem{(i)} $E_s$, if not empy, is a direct sum of stable
bundles, each of slope $\tau\cdot {\hbox{\smalltype Vol}(\Sigma)\over
4\pi}$;
\itemitem{(ii)} $E_\phi$ contains the section $\phi$ and
$(E_\phi,\phi)$ is $\tau$-stable, where $E_\phi$ has the
holomorphic structure induced from $\dbare$.

Notice that the split case $E=E_\phi\oplus E_s$
cannot occur unless $\tau\cdot
{\hbox{\smalltype Vol}(\Sigma)\over
4\pi}$ corresponds to the slope of a subbundle, i.e. unless $\tau\cdot
{\hbox{\smalltype Vol}(\Sigma)\over
4\pi}$\ is a rational number with
denominator less than the rank of $E$.  Hence,
for generic values of $\tau$ the summand $E_s$ is empty, and the
$\tau$-stable pairs comprise the set
$$\Vt=\Bigl\{ (\dbf)\in\H :
\Lambda\curvdh -{ \sqrt{-1}\over 2}\phi
\otimes\phi^\ast=-\sqrt{-1}{\tau\over 2}{\bf I}\ \  \hbox{for some metric } H
\Bigr\}\;.\leqno(2.4)$$

An important feature of the vortex equation is its interpretation as a
symplectic
moment map condition.  This comes about as follows.
If a Hermitian bundle metric, $H$
say, is fixed on $E$, then $\cal H$
acquires a natural symplectic
structure, i.e. the one coming from the usual symplectic structures on ${\cal
C}$ and $\Omega^0(E)$\ (cf. [B-D1]).  Moreover,
the unitary gauge group $\G$ acts {\it
symplectically}  and the  moment map for this action is exactly the left hand
side
of the vortex equation, viz.
$$\mmap (\dbf)=\Lambda\curvdh - {\sqrt{-1}\over
2}\phi\otimes\phi^\ast\;.\leqno(2.5)$$
It follows that for generic $\tau$,
the space $\Vt$ is the saturation of
$\Psi^{-1}(-{{\sqrt{-1}\tau}\over 2}\bf I)$, i.e. consists of all the
$\Gc$-orbits through the holomorphic pairs in $\Psi^{-1}(-{{\sqrt{-1}\tau}\over
2}\bf I)$.  We thus get two descriptions of the moduli spaces
$$\Btau = \Vt/\Gc=\Psi^{-1}\left(-{{\sqrt{-1}\tau}\over 2}\bf I\right
)/\G
\leqno(2.6)$$
of $\tau$-stable pairs.
The first description gives the
complex structure and the second gives the
symplectic structure.

Our new space $\Bhat$
will similarly have two descriptions; one as a complex
orbit space, and one as a
symplectic reduction from a moment map. To see what  the
moment map should be we start from the fact that the space
is designed to contain all holomorphic pairs that are
$\tau$-stable for
\it some \rm value of $\tau$.  Hence, using the vortex
equation
characterization of $\tau$-stability,
it is evident that the pairs in $\Bhat$ are
characterized by the condition that for some metric $H$,
$$\Lambda\curvdh -{ \sqrt{-1}\over 2}\phi
\otimes\phi^\ast=\hbox{const.}\,{\bf I}.$$
If we let  $\kform 0 {\End E}_0$
denote the $L^2$-orthogonal complement of the
constant multiples of the identity in $\kform 0 {\End E}$, and let
$\pi^{\perp}:\kform 0 {\End E}\to\kform 0{\End E}_0$\ denote
the orthogonal projection, then the defining condition for $\Bhat$
becomes
$$\pi^{\perp}(\Lambda\curvdh -{ \sqrt{-1}\over 2}\phi
\otimes\phi^\ast)=\pi^{\perp}\mmap(\dbf)=0.\leqno(2.7)$$
This can be realized as a moment map condition
if we replace the full unitary
gauge  group $\G$ by a  subgroup $\G_0\subset\G$
which has a $U(1)$ quotient
and whose Lie algebra is the
ortho-complement (with respect
to the $L^2$-metric) of the constant multiples of
the identity. Denoting the
new moment map by $\Psi_0$, the ``master space"
$\Bhat$ will
then correspond to the reduced
zero set $\Psi_0^{-1}(0)/\G_0$.  Finally, to  obtain
the complex description of the space we must now
consider the saturation of
$\Psi_0^{-1}(0)$, not with respect to
orbits of the full complex gauge group,
but with respect to orbits of the
subgroup corresponding to
the complexification of $\G_0$.

\centerline{\it \S 2.2 The subgroups of $\G$\ and $\Gc$}

A key element in the construction of $\Bhat$ is thus
the choice of an appropriate subgroup of
$\G$\ (or $\Gc$).  At least in the connected
component of the identity,
the discussion above makes it clear that the
subgroup we want is the
unique connected subgroup corresponding to the Lie
subalgebra $\Omega^0(\End E)_0$.
In the image of the exponential map, this subgroup can also be
described as the kernel of the homomorphism
$$\chi(\exp (u)) = \exp \left(\int_X \Tr(u)\right)\;.$$
In order to define the full subgroup we need to extend this
homomorphism to all of $\G$.  The technical details are as follows,
beginning with

\proclaim Lemma 2.2.  Let $f:\Sigma\to {\Bbb C}^\ast$
be a smooth map in the
connected component of the identity in $\hbox{Map}(\Sigma,
{\Bbb C}^\ast)$.
Then there exists a unique $\chi_1(f)\in{\Bbb C}^\ast$ and
$u:\Sigma\to{\Bbb C}$ satisfying $f=\chi_1(f)\exp u$ and $\int_\Sigma
u=0$.  Moreover, if $f'$ is another such map,
$\chi_1(ff')=\chi_1(f)\chi_1(f')$.

\noindent {\it Proof.}
Let us first prove uniqueness.  If
$$f=\chi_1(f)\exp u_1 =\chi_2(f)\exp u_2\;,$$
then $u_1-u_2$ must be constant.  But then $\int_\Sigma u_1-u_2=0$
implies the constant is zero, and so $u_1=u_2$ and
$\chi_1(f)=\chi_2(f)$. Similarly, $\chi_1(ff')=\chi_1(f)\chi_1(f')$.
To prove existence, let $\Sigma={\Bbb H}/\Gamma$ where ${\Bbb H}$ is
the upper half plane in $\Bbb C$ and $\Gamma\subset\hbox{PSL}(2,{\Bbb
R})$ is the uniformizing group.  Then $f$ lifts to $\tilde f:{\Bbb
H}\to {\Bbb C}^\ast$ satisfying $\tilde f(\gamma z)=\tilde f(z)$ for
all $\gamma\in \Gamma$.  Choose a point $z_0\in {\Bbb H}$, and let
$u_0:{\Bbb H}\to {\Bbb C}$ be defined by
$$u_0(z)=\int_{z_0}^z{d\tilde f\over \tilde f}\;.\leqno(2.8)$$
Since $d\tilde f/\tilde f$ is closed, this is independent of the
path from $z_0$ to $z$.  Moreover, since
$\Gamma\subset\hbox{PSL}(2,{\Bbb R})$ we have
$$\int_{\gamma z_0}^{\gamma z}{d\tilde f\over \tilde f}=
\int_{z_0}^z{d\tilde f\over \tilde f}\;,$$
from which it follows that
$$u_0(\gamma z)=u_0(z)+
\int_{z_0}^{\gamma z_0}{d\tilde f\over \tilde f}\;,$$
for all $z\in {\Bbb H}$, $\gamma\in\Gamma$.  The second term is just
the winding number of $f$ about the cycle defined by $\gamma$, and
this vanishes since $f$ is assumed to be connected to the identity.
Thus $u_0$ descends to a map $u_0:\Sigma\to{\Bbb C}$, and clearly
$f=\hbox{const.}\,\exp u_0$. Normalizing
$u=u_0-{1\over\hbox{\smalltype Vol}(\Sigma)}\int_\Sigma u_0$, we
obtain $\chi_1(f)$.  It is now easily verified that both $u$ and
$\chi_1(f)$ are independent of the choice of point $z_0$.

Let $\Gci\subset \Gc$ denote the connected component of the
identity, and let $\Upsilon$ denote the quotient group of
components.  Then $\Upsilon$ is a free abelian group on $2g$
generators corresponding to $H_1(\Sigma, {\Bbb Z})$ (see [A-B], p.
542).  We can find a splitting of the exact sequence
$$1\lra\Gci\lra\Gc\lra\Upsilon\lra 1\leqno(2.9)$$
which realizes $\Gc$ as a direct product
$$\Gc\simeq\Gci\times \Upsilon\;,\leqno(2.10)$$
where the isomorphism is given by $(g,h)\mapsto gh$.  Using Lemma
2.2 and the isomorphism (2.10), we can now define a character on
$\Gc$ as follows:  First, for $g\in\Gci$, $\det
g:\Sigma\to\CBbb^\ast$ is in the connected component of the
identity, and so we may set $\chi(g)=\chi_1(\det g)$.  Then we
extend $\chi$ to $\Gci\times \Upsilon$ by $\chi(g,h)=\chi(g)$.
This defines a homomorphism $\Gc\to\CBbb^\ast$.

\proclaim Definition 2.3.  Let $\Gco$ be the kernel of the
character $\chi:\Gc\to{\Bbb C}^*$ defined as above.
Let $\G_0\subset\G$\ be defined by
$\G_0=\Gco\cap\G$.

Note that a different choice of splitting, or isomorphism (2.10), will
give rise to an isomorphic group $\tilde\Gco$ with the same
connected component of the identity as $\Gco$.  The following is
immediate from the definition:
\proclaim Proposition 2.4.
The groups $\G_0$\ and $\Gco$ have the structure
of  Fr\' echet  Lie groups with Lie algebras
$$\hbox{Lie}\,\Gco=\Omega^0(\End E)_0=\{ u\in\Omega^0(\End E) :
\int_\Sigma \Tr u=0\}\;,\leqno(2.11)$$
and
$$\hbox{Lie}\,\G_0=\Omega^0(\hbox{ad}\, E)_0=\{ u\in\Omega^0(
\hbox{ad}\, E) :
\int_\Sigma \Tr u=0\}\;.\leqno(2.12)$$

\medskip
\centerline{ \it \S 2.3 Local complex structure}

We now proceed with the construction of the master space
and begin with the
construction of $\Bhat$ as a complex manifold.
As a complex manifold, $\Bhat$
is essentially the orbit space for the $\Gco$
action on $\cal H$.
The construction is thus a relatively small modification
of the procedures used for
the orbit space $\cal H/\Gc$.  In that case, the
obstructions to having
a smooth manifold structure, as well as the description
of the tangent spaces,
come from the cohomology of the deformation  complex
$\Cx$:
$$0\mapright{} \kform 0 {\End E}\mapright{d_1} \pqform 0 1 {\End
E}\oplus
\kform 0 E \mapright{d_2} \pqform 0 1 E\mapright{} 0\;,\leqno(2.13)$$
where the maps $d_1$, $d_2$ are given by
$$\eqalign{d_1(u)&= (-\dbare u, u\phi)\cr
	d_2(\alpha, \eta) &= \dbare\eta +
	\alpha\phi\;.\cr}\leqno(2.14)$$
The salient features of this complex are contained in the following

\proclaim Proposition 2.5.  Let $(\dbf)$ be an element of $\H$.  Then
\itemitem{(i)}  $\Cx$ is an elliptic complex;
\itemitem{(ii)} If $d>R(2g-2)$, then $H^2(\Cx)=0$ whenever
$\phi\neq 0$ or
$E^{\dbare}$ is semistable;
\itemitem{(iii)}  The Euler characteristic of the complex $\Cx$ is
given by
$$\chi(\Cx)=\chi(\End E)-\chi(E)\;.$$

\noindent{\it Proof.}
This is proven in [B-D1], except there the case $\phi=0$\
is excluded and  an extra assumption
is made which ensures that $E^{\dbare}$
is always semistable.
It is clear that, under the assumption that $d>R(2g-2)$,
the case $\phi=0$
can be included without any modification to the proof.
Furthermore, as shown in [Th],
when $\phi\ne 0$, the vanishing of $H^2(\Cx)=0$
follows in general since the map
$$H^0(K\otimes E^\ast)\mapright{\otimes\phi}H^0(K\otimes\End E)$$
is injective.  By Serre duality this is equivalent to the surjectivity of the
map\break
$H^1(\End E)\to H^1(E)$\ in the long exact sequence
$$0\mapright{}
H^0\mapright{}H^0(\End E)\mapright{}H^0(E)\mapright{}H^1
\mapright{}H^1(\End E)\mapright{}H^1(E)\mapright{} H^2\mapright{} 0\;,
\leqno(2.15)$$
where $H^i=H^i(\Cx)$.
The proofs of (i) and (iii) are as in [B-D1].

For the construction of $\Bhat$, we need to restrict to
the following subcomplex of  $\Cx$:
$$ 0\mapright{} \kform 0 {\End E}_0
\mapright{d_1} \pqform 0 1 {\End
E}\oplus \kform 0 E
\mapright{d_2} \pqform 0 1 E\mapright{} 0\;,\leqno(2.16)$$
which we will refer to as $\Cxo$.
It is worth pointing out that
the adjoint $d_1^{*,0}$\ in $\Cxo$\ is related to
the $d_1^*$\ in $\Cx$ by
$$d_1^{*,0}= \pi^{\perp} d_1^*,\leqno(2.17)$$
where $\pi^{\perp}$
gives the orthogonal projection onto $\kform 0 {\End
E}_0$.
This affects the determination of the harmonic  1-cocycles in the two
complexes.  In fact, we have
\proclaim Proposition 2.6.  Let $(\dbf)$ be an element of $\H$.  Then:
\itemitem{(i)}  $\Cxo$ is a Fredholm complex;
\itemitem{(ii)} $H^2(\Cxo)= H^2(\Cx)$;
\itemitem{(iii)} Either
$$\left\{   \eqalign{ H^1(\Cxo)&\simeq H^1(\Cx)\oplus {\Bbb C}\cr
H^0(\Cxo)&= H^0(\Cx)\cr}\right.$$
or,
$$\left\{  \eqalign {H^1(\Cxo)&= H^1(\Cx)\cr
H^0(\Cxo)\oplus{\Bbb C} &\simeq H^0(\Cx)\cr}\right.$$
\itemitem{(iv)} $\chi(\Cxo)=\chi(\End E)-\chi(E)-1\;.$

\noindent{\it Proof.}
(i), (ii) follow immediately from the definition of $\Cxo$.\hfil\break
(iii)  $H^1(\Cx)$\ and $H^1(\Cxo)$
are related by the short exact sequence
$$ 0\to H^1(\Cx)\to H^1(\Cxo)\to\Bbb C\to 0.$$
Similarly, the zero-th cohomology groups are related by
$$ 0\to H^0(\Cxo)\to H^0(\Cx)\to\Bbb C\to 0,$$
where the map is orthogonal projection in $\kform 0 {\End E}$.
Here by $\Bbb C$
we mean the constant multiples of the identity in $\kform
0{\End E}$. The desired conclusion now follows from the fact that
the map $d_1^\ast: H^1(\Cxo)\to{\Bbb C}$ is surjective if and only
if $\pi : H^0(\Cx)\to\Bbb C$ is zero.
(iv) follows from (ii) and (iii).

\proclaim Corollary 2.7.  Let $\H^\ast\subset \H$ denote the
subspace of all $(\dbf)\in\H$ such that $E^{\dbare}$ is semistable if
$\phi=0$.  Then $\H^\ast$ is a smooth submanifold of $\aspace$.

\noindent{\it Proof.}  Consider the map
$F:\aspace\longrightarrow\kform 1 E$ defined by
$F(\dbare,\phi)=\dbare(\phi)$.  The derivative of $F$ at $(\dbare,
\phi)$ is given by
 $$\delta F_{\dbf}(\alpha,\eta)=\dbare
 \eta+\alpha\phi=d_2(\alpha,\eta)\;.\leqno(2.18)$$
Let ${\cal C}_{ss}$ denote
the semistable holomorphic structures on $E$.  Then
provided that $(\dbf)$ does not belong to the closed subspace
$({\cal C}-\ss)\times \{0\}\subset \aspace$, it follows from Proposition 2.5
(ii) that $\delta F_{\dbf}$ is onto.
Hence by the Inverse Function Theorem

$$\H^\ast= F^{-1}(0)\cap\left\{ \aspace-({\cal
C}-\ss)\times\{ 0\}\right\}\leqno(2.19)$$
is a smooth submanifold of $\aspace$.

\proclaim Definition 2.8.  A pair $(\dbf)\in\H^\ast$ is called
simple if $H^0(\Cxo)=0$.  Let $\H_\sigma$  denote the subspace of
simple pairs in $\H^\ast$.

Clearly, $\H_\sigma$ is an
open subset in $\H^\ast$ and is therefore a submanifold.
Now by identifying
$H^1(\Cxo)$ with the tangent space to the space of
orbits of $\Gco$, we have the following theorem whose proof
is in all essentials the same as the proof of the analogous result
in Section 3 of [B-D1]:

\proclaim Theorem 2.9.  $\H_\sigma/\Gco$ is a complex V-manifold
(possibly non-Hausdorff) of complex dimension $d+1+(R^2-R)(g-1)$.
Moreover, we have the identification
$$T_{[\dbf]}(\H_\sigma/\Gco)\simeq H^1(\Cxo)\;.\leqno(2.20)$$

The space $\H_\sigma/\Gco$
is almost, but not quite, the master space $\Bhat$.
The problem is the possibly non-Hausdorff
nature of the space.  This can be
traced back (as in all such moduli space problems)
to the fact that the
definition of a ``simple" pair is too weak; it
needs to be replaced by a
concept of ``stability",
i.e. we need to restrict from $\H_{\sigma}$\ to some
suitable analogue of the $\Vt$
used in the construction of $\Btau$.  We recall
briefly the
procedure for the $\Btau$.  For generic values of $\tau$, i.e.
$\tau$ not equal to a rational number with
denominator
less than $R$, the set of $\tau$-stable pairs in $\cal H$\ is given
precisely by the set $\Vt$ defined in (2.4).
Furthermore for such $\tau$, the
set $\Vt$
is open subset in,  and therefore a submanifold of,
$\cal H^\ast$.  The
quotient
$\Vt/\Gc$ is homeomorphic to
$\mmap^{-1}( -{\sqrt{-1}\over 2}\tau\cdot\bf I)/\G$, which is a
Hausdorff,
smooth, symplectic manifold.  In [B-D1] these last assertions are
proven in the case where the degree of $E$
is large and $\tau$ is small
(Assumptions (1) and (2)).  However,
in view of Corollary 2.7
it is  clear that the conclusion holds more
generally.  That is,

\proclaim Proposition 2.10.
For any degree and any generic value of $\tau$,
the space $\Vt/\Gc\simeq\mmap^{-1}(-{\sqrt{-1}\over 2}
\tau\cdot\bf I)/\G$, if nonempty, is a
smooth compact K\"ahler manifold. It is $\Btau$,
the moduli space of $\tau$-stable
pairs.

\noindent{\it
Proof}. The purpose of the two assumptions in [B-D1] was primarily to
ensure that the subset of $\cal H$\ used in the
construction of the moduli
spaces
was a submanifold of $\aspace$.  However,
by Corollary 2.7 we see that this can
be
achieved without any restriction on
$\tau$ by using the space $\cal H^\ast$.
The
constraint on the degree of the bundle can also be relaxed; for the
construction of
 the moduli spaces of $\tau$-stable pairs all that is required is the
vanishing of
 $H^2(\Cx)$\ when $\phi\ne 0$,
 or equivalently the surjectivity of the map
$$\delta F_{\dbf}:\pqform 0 1 {\End E}\longrightarrow\kform 1 E\;\leqno(2.21)$$
when $\phi\ne 0$.

It will also be important to know what the allowed range for $\tau$
is.  By
taking the trace of the
vortex equation and integrating over the base manifold
$\Sigma$, one sees that there can be no solutions unless
$$\tauhat\ge{d\over R}\quad.\leqno(2.22)$$
Furthermore, an
upper bound on $\tau$\ can be obtained by looking at the quotient
$E/[\phi]$, where $[\phi]$\ is the line subbundle
generated by the section
$\phi$.  If the $\tau$-vortex equation is satisfied, then (cf. [B1])
$$\tauhat\le{\hbox{deg}\,(E/[\phi])\over \rank(E/[\phi])}\quad ,$$
and hence, since $\hbox{deg}\,([\phi])\ge 0$,
$$\tauhat\le{d\over R-1}\quad .\leqno(2.23)$$
In fact,
solutions corresponding to all intermediary values of $\tau$
between these extreme points
can be constructed
(see  Propositions  3.2 and 2.18), and  we thus have

\proclaim Proposition 2.11.  There is a solution
to the $\tau$-Vortex equation
if and only if $\tauhat$ is in the closed interval $[d/R , d/
R-1]$.

\medskip
\centerline{\it\S 2.4 Symplectic structure and global complex structure}

We now consider the symplectic (K\"ahler) structure on
$\H_\sigma/\Gco$.  In particular we need the $\G_0$-moment map, the symplectic
reduction of its zero set, and the $\Gco$-saturation of this zero set.

\proclaim Proposition 2.12.  A moment map for the action of $\G_0$
on $\H^\ast$ is given by
$$\mmapo(\dbf)=\pi^{\perp}\mmap(\dbf)=\mmap(\dbf)-{1\over
R\cdot\hbox{\smalltype
Vol}(\Sigma)}\int_\Sigma \Tr\mmap(\dbf)\,\cdot{\bf I}\;.\leqno(2.24)$$

\noindent{\it Proof.}  Let $\jmath : \G_0\to\G$ denote the
inclusion.  Then a moment map for $\G_0$ is given by
$\jmath^\ast\Psi$, where $\jmath^\ast :
(\Lie\G)^\ast\to (\Lie\G_0)^\ast$.
Using the $L^2$-inner
product on $\kform 0 {\End E}$ to identify the
Lie algebras with their duals, we obtain $\Psi_0$.

\proclaim Definition 2.13.  Let
$$\Bhat=\mmapo^{-1}(0)/\G_0\leqno(2.25)$$
 denote
the Marsden-Weinstein
reduction of $\H^\ast$ by the symplectic action of $\G_0$.
In addition, let $\H_0\subset
\H^\ast$ denote the subspace of $\H^\ast$ where $\G_0$ acts with at
most finite stabilizer.  We then define
$$\Bhat_0= \mmapo^{-1}(0)\cap \H_0/\G_0\;.\leqno(2.26)$$

\proclaim Proposition 2.13.
The quotient $\Bhat$ is a compact, Hausdorff
topological space.  The quotient $\Bhat_0$ is a Hausdorff symplectic
V-manifold.

\noindent {\it Proof.}
To prove the compactness of $\Bhat$, let $(D_i,\phi_i)$ be a
sequence in $\Psi_0^{-1}(0)$. We will have to find $g_i\in\G_0$ such
that $(g_i(D_i), g_i\phi_i)\to (D,\phi)\in \Psi_0^{-1}(0)$.  Let us
write
$$\Lambda F_{D_i}+{\sqrt{-1}\over 2}\phi_i\otimes
\phi_i^\ast=-{\sqrt{-1}\over 2}\tau_i$$
for real numbers $\tau_i\in [d/R, d/R-1]$.  By passing to a
subsequence we may assume that $\tau_i\to\tau$.  The rest of the
compactness argument follows as in [B-D1], Proposition 5.1, and the
observation that $g_i$ can be chosen to be in $\G_0$, since the
constant central elements of the gauge group act trivially on $\cal
C$.  The Hausdorff property follows exactly the same way as in
[B-D1], Proposition 5.4.  Finally, for the symplectic structure we
use the symplectic reduction theorem for the group $\G_0$ (cf.
[B-D1], Theorem 4.5).  The only difference is that one must replace
the complex $M_\phi^{\dbare}$ of [B-D1] by the subcomplex
$$
0\mapright{} \kform 0 {\ad E}_0\mapright{D_1} \pqform 0 1 {\End
E}\oplus \kform 0 E \mapright{} \kform 0 {\ad E}\oplus
\pqform 0 1 E\mapright{} 0\;.$$
This completes the proof.

To complete the construction
of the master space as a complex manifold we
make the following
\proclaim Definition 2.15.
Let ${\cal V}_0\subset\cal H$\ denote the subset of
$\Gco$-orbits through points in $\Psi_0^{-1}(0)$, i.e.
$${\cal V}_0
= \{(\dbf)| \Psi_0( g(\dbare,\phi))=0   \hbox{ for some } g \in
\Gco\}.$$

It is easily seen that ${\cal V}_0\cap{\cal H}_\sigma$
is an open subset of ${\cal
H}_{\sigma}$, and thus is a submanifold.
Also, using the same techniques as those applied to ${\cal V}_\tau$
in [B-D2] it can be shown that ${\cal V}_0$ and ${\cal V}_0\cap
{\cal H}_\sigma$ are
connected.  Indeed, the only essential difference in the argument is
that now the projection to the set of holomorphic structures
 includes unstable structures.  Since these add
sets of small codimension they do not affect the connectedness.
Finally, there is clearly a  bijective
correspondence between $\Bhat_0$\ and ${\cal V}_0\cap{\cal H}_\sigma
/\Gco$.  Combining
Theorems 2.9 and 2.13 we thus obtain
\proclaim Theorem 2.16.
$\Bhat_0={\cal V}_0\cap{\cal H}_\sigma/\Gco$ is a smooth, Hausdorff,
K\"ahler manifold
of dimension \hfil\break
$d+1+(R^2-R)(g-1)$.

\centerline{\it
\S 2.5  $S^1$-action , Morse function, and reduced level sets}

The most important
feature of the master space $\Bhat$ is the fact
that it carries an $S^1$ action.  This comes from the quotient
$\U\simeq\G/\G_0$, and the action on $\Bhat$ is given by
$$e^{i\theta}\cdot[\dbf]=[ \dbare, g_\theta\phi]\;.\leqno(2.27)$$
Here $g_\theta$ denotes the gauge transformation
$\hbox{diag}(e^{i\theta/R},\ldots,e^{i\theta/R})$.
Notice that $g_\theta$ itself depends on the choice of an $R$-th
root of unity but that the action is well-defined and independent of
this choice, since if $h=
\hbox{diag}(e^{2\pi i/R},\ldots,e^{2\pi i/R})$, then $h\in \G_0$ and
$$[\dbare, h\phi]=[h^{-1}(\dbare),\phi]=[\dbf]\;.$$
The
following summarizes the basic properties of this action:

\proclaim Proposition 2.17.  (i) The action of U(1) on $\Bhat_0$ is
holomorphic and symplectic; the moment map for the action is given
by
$$\hat f[\dbf]=-2\pi i\left({\Vert \phi\Vert^2\over 4\pi
R}+\mu(E)\right)\;.\leqno(2.28)$$
(ii) The action extends continuously to $\Bhat$  as does the moment
map $\hat f$.

\noindent {\it Proof.}  The only statement that needs to be verified
is the computation of the moment map.  First, observe that
there is no natural splitting of the exact sequence of groups
$$1\mapright{}\G_0\mapright{}\G\mapright{\chi}\U\mapright{} 1\;.$$
However, on the level of Lie algebras there is a splitting $\lambda
: \hbox{Lie}\,\U\to\hbox{Lie}\,\G$ given by $\lambda{(\xi)}=\xi/R$,
where $\xi\in \hbox{Lie}\,\U$ is identified with the constant
infinitesimal gauge transformations.  Under the identification with
dual Lie algebras given by the $L^2$ metric,
$$\lambda^\ast : (\hbox{Lie}\,\G)^\ast\lra
(\hbox{Lie}\,\U)^\ast $$
is given by
$$\lambda^\ast(g)=\int_\Sigma\Tr g\;.$$
Now given $[\dbf]\in\Bhat_0$, choose
a representative $(\dbf)$ satisfying $\mmapo(\dbf)=0$.  It is
straightforward to compute that
$$\hat f[\dbf]={1\over R}\,\lambda^\ast\circ\mmapo(\dbf)=-2\pi
i\left({\Vert\phi\Vert^2\over 4\pi R}+\mu(E)\right)\;.\leqno(2.29)$$
Notice that if we represent a point $[\dbf]$\ in $\Bhat_0$\ by
$(\dbf)\in\mmapo^{-1}(0)$, then
$$ {-\hat f[\dbf]\over 2\pi i}=\tauhat\leqno(2.30)$$
if and only if
$$\mmap(\dbf)={-\sqrt{-1}{\tau}\over{2}}\bf I\ \rm.\leqno(2.31)$$

\noindent For convenience we let $f:\Bhat\to{\Bbb R}$ denote the
function
$$f=-{1\over 2\pi i}\hat f\;.\leqno(2.32)$$
We now have
\proclaim Proposition 2.18.  (i) The image of $f$ is the interval
$[d/R,d/R-1]$.\hfil\break
(ii) The critical points of $f$ on $\Bhat_0$ are
precisely the fixed points of the U(1) action.  The critical values
of $f$ coincide with the image under $f$ of the fixed point set of
the U(1) action on $\Bhat$.\hfil\break
(iii)  Let
$\hat\tau=\tau\cdot{\hbox{\smalltype Vol}(\Sigma)\over 4\pi}$
be a regular value of $f$.  Then
the reduced space $f^{-1}(\hat\tau)/\U$ is a K\"ahler manifold which
can be identified with the moduli space of $\tau$-vortices in degree
$d$ and rank $R$, i.e.
$$f^{-1}(\hat \tau)/\U =\Btau\;.$$

\noindent {\it Proof.} (i)  By the
comment after Proposition 2.17, $\tau$ is in
the image of $f$ if and only if the equation
$\mmap(\dbf)={-\sqrt{-1}{\tau}\over{2}}\ \bf I\rm$
has a solution.
Hence, by Theorem 2.1 the range for $\tau$ is in $[d/R,d/R-1]$.
Now the endpoints of this interval are included in the image, since
explicit elements of the pre-image can be constructed (see
Proposition 3.2).  The result then follows from
the connectedness of ${\cal V}_0\cap {\cal H}_\sigma$
(see the remark following Definition 2.15).
(ii)  This follows from the fact that $\hat f$\ is a moment map for the $U(1)$\
action.
(iii) Given $[\dbf]\in f^{-1}(\hat\tau)/\U$, choose a
representative $(\dbf)\in\Psi_0^{-1}(0)$.  Then by definition of the
moment maps we obtain the pair of equations
$$\eqalign{\Lambda\curvdh-{\sqrt{-1}\over 2}\phi\otimes\phi^\ast&=
{1\over R}\int_\Sigma\Tr\left(
\Lambda\curvdh-{\sqrt{-1}\over 2}\phi\otimes\phi^\ast\right)\cdot{\bf
I}\cr
{\Vert\phi\Vert^2\over 4\pi R}+\mu(E)&=\hat\tau\;,\cr}\leqno(2.33)$$
which are clearly equivalent to the $\tau$-vortex equation (2.3).
Thus, taking the class of $(\dbf)$ in $\Btau$ defines a map
$$f^{-1}(\hat\tau)/\U\lra \Btau\;.\leqno(2.34)$$
The inverse map, i.e. mapping the class of $(\dbf)$\ in
$\Btau=\mmap^{-1}({-\sqrt{-1}{\tau}\over{2}}\ \bf I\rm)/\G$\ to the class of
$(\dbf)$\ in $f^{-1}(\hat\tau)/U(1)$, is well defined, and hence the map is a
bijection.
It is easily checked that for non-critical values of $\tau$ this
map is indeed an isomorphism of K\"ahler manifolds.

Let ${\cal T}$ denote the set of critical values of $f:\Bhat_0\to
{\Bbb R}$.  We shall see in the next section that ${\cal T}$
consists of the rational numbers in $[d/R, d/R-1]$ which appear as
slopes of subbundles of $E$.  Implicit in Proposition 2.18 is the
statement that $\Bhat\setminus f^{-1}({\cal T})=\Bhat_0\setminus
f^{-1}({\cal T})$ is smooth.  This is clear, since a point
$(\dbf)\in \mmapo^{-1}(0)$ can have non-trivial isotropy in $\G_0$
only if the bundle $E$ splits holomorphically or if $\phi\equiv 0$.
In rank two we can say more:

\proclaim Proposition 2.19.  Suppose $R=2$.  Then
$$\Bhat\setminus f^{-1}(d/2)=\Bhat_0\setminus f^{-1}(d/2)\;,\leqno(2.35)$$
and this space is smooth.  If the degree $d$ is odd, then
$\Bhat=\Bhat_0$, and this space has at most ${\Bbb Z}_2$ quotient
singularities along $f^{-1}(d/2)$.

\noindent {\it Proof.}  Suppose $(\dbf)\in \mmapo^{-1}(0)$.  Then
$(\dbf)\in\mmap^{-1}(-\sqrt{-1}{\tau\over 2}
\bf I\rm)$ for some $\tau$.  By
Theorem 2.1,
$\phi\equiv 0$ if and only if $\hat\tau=d/2$.  Now suppose that
$g\in \G_0$, $g\neq I$, and $g\cdot(\dbf)=(\dbf)$.  If
$\phi\not\equiv 0$, then $E$ must split holomorphically as
$E=E_\phi\oplus E_s$ with $\phi\in H^0(E_\phi)$ and $g=(1,\tilde
g)$.  But $\rank(E_s)=1$ implies $\tilde g$ is constant, and since
$g\in\G_0$ we must have $\det(g)=\tilde g=1$.  This proves that
the stabilizer for points away from $\Psi_{d/2}^{-1}(0)$ is trivial,
and therefore
$$\Bhat\setminus f^{-1}(d/2)=\Bhat_0\setminus f^{-1}(d/2)$$
is smooth.  If in addition we assume that $d$ is odd, then
$(\dbf)\in \Phi_{d/2}^{-1}(0)$ implies that $\phi\equiv 0$ and $E$
is stable; hence the stabilizer consists of $\pm I$.  This completes
the proof.

\beginsection 3.  Critical  sets

To further our understanding of  the master space $\Bhat$ we now
describe the level sets $f^{-1}(\hat\tau)$ for the critical values
$\hat\tau$ in the interval $[d/R, d/R-1]$.

\proclaim Definition 3.1.  Let $\hbox{Fix}(\Bhat)$ denote the U(1)
fixed point set in $\Bhat$.  For a critical value
$\hat\tau=\tau\cdot{\hbox{\smalltype Vol}(\Sigma)\over 4\pi}$, let
$$\Ztau=f^{-1}(\hat\tau)\cap\hbox{Fix}(\Bhat)\;.\leqno(3.1)$$

\proclaim Proposition 3.2.  The level set corresponding to the
minimum is precisely the moduli space of semistable bundles of
degree $d$ and rank $R$, i.e.
$$f^{-1}(d/R)={\cal M}(R,d)\;.\leqno(3.2)$$
The level set corresponding to the maximum is precisely the moduli
space of semi-stable bundles of degree $d$ and rank $R-1$, i.e.
$$f^{-1}(d/R-1)={\cal M}(R-1,d)\;.\leqno(3.3)$$

\noindent {\it Proof.} By Proposition 2.18, $f^{-1}(d/R)$ consists
precisely of $\tau$-semistable pairs where $
\phi\equiv 0$.  But then $\tau$-semistable is equivalent to
semistable.
For $f^{-1}(d/R-1)$, suppose $(\dbf)$\ supports a solution to the
$\tau$-vortex equation with
$\tau=d/R-1$.  By Theorem 2.1, $(\dbf)$ is
either a stable pair,
or splits holomorphically.  The former
is not possible,
as can be seen by applying the second condition of Definition
1.1 to $[\phi]$, the line subbundle generated by $\phi$.  This yields
$\hbox{deg}\,
([\phi])<0$.  In fact the only possibility is that $
\hbox{deg}\,([\phi])=0$ and
thus $\mu(E/[\phi])=\tau$.  It follows (cf. Theorem 2.1 and
Section 2.4 of [B1])
that $\phi$\ is a constant section of a trivial line subbundle and $E={\cal
O}\oplus
E_s$, where $E_s$ is a semistable bundle of degree $d$, rank $R-1$. In fact,
the summand $E_s$ is a direct sum of stable bundles all of the same slope,
i.e. $E_s$\ is equal to the graded
bundle in its S-equivalence in ${\cal M}(R-1,
d)$.  As above, one can check that the map ${\cal O}\oplus
E_s\mapsto E_s$\ does indeed give a bijective correspondence between
$f^{-1}(d/(R-1))$\ and ${\cal M}(R-1,d)$.

Next, we describe the level sets for the intermediate values of
$\hat\tau$.  Let $\hat\tau=p/q\in(d/R, d/(R-1))$ be a critical value.
Then any pair $(\dbf)\in\Ztau$ has $\muM=p/q=\mum$, and the
bundle splits holomorphically as $E^{\dbare}=E_\phi\oplus E_{ss}$,
where
\itemitem{(i)} $\phi\in H^0(E_\phi)$,
\itemitem{(ii)} $(E_\phi,\phi)$ is a $\tau$-stable pair,
\itemitem{(iii)} $E_{ss}$ is a direct sum $\bigoplus_i E_i$ of stable
bundles, all of slope $\hat\tau$.

\proclaim Lemma 3.3.  Fix a critical value $p/q\in (d/R, d/(R-1))$.
Suppose that $E^{\dbare}=E_\phi\oplus E_{ss}$ is part of a
non-$\tau$-stable pair in $f^{-1}(p/q)$ as above.  Let the degree
and rank of $E_\phi$ be $(R_\phi, d_\phi)$ and those of $E_i$ be
$(R_i, d_i)$.  Then we have the following constraints:
\itemitem{(i)} $d_i/R_i =(d-d_\phi)/(R-R_\phi)=p/q$;
\itemitem{(ii)} $\sum_i R_i = R-R_\phi$;
\itemitem{(iii)} $d_\phi/R_\phi < p/q < d_\phi/(R_\phi-1)\;.$

\noindent {\sl Conversely, given any stable pair $(E_\phi, \phi)$ and
set of stable bundles $E_i$ such that the conditions above are
satisfied, we obtain a representative for a
fixed point $(E_\phi\oplus\bigoplus_i E_i,
\phi)  $ in the critical level set corresponding to $p/q$.}

We see that to a fixed point $(\dbf)\in\Ztau$ we can assign an
$(n+2)$-tuple of integers $\tuple$.  The
degrees $d_i$ are then determined by condition (i) of Lemma 3.3.
For $\hat\tau=p/q\in(d/R, d/(R-1))$, let
$$\eqalign{I_\tau=\Bigl\{\tuple
&\in \ZBbb^{n+2} :
 \hbox{ with }
{p\over q}=\hat\tau={d-d_\phi\over R-R_\phi}\ \cr
& \hbox{ conditions (ii) and (iii) of Lemma 3.3 are satisfied}
\;\Bigr\}\;,\cr}\leqno(3.4)$$
and set
$${\cal Z}\tuple=\left\{ (\dbf)\in\Ztau :
E^{\dbare}=E_\phi\oplus\bigoplus_{i=1}^n E_i\right\}\leqno(3.5)$$
where $\degree(E_\phi)=d_\phi$, $\rank(E_\phi)=R_\phi$, and
$\rank(E_i)=R_i$ for $i=1,\ldots, n$.
We can extend this notation to include
the extreme values $\hat\tau= d/R$ and
$\hat\tau= d/(R-1)$. This requires the convention that
\itemitem{(a)} if $\hat\tau= d/R$, then
$d_{\phi}=R_{\phi}=0$\ and $\sum_i R_i = R$,
\itemitem{(b)}  if $\hat\tau= d/(R-1)$, then
$d_{\phi}= 0$, $R_{\phi}=1$
and conditions (ii) and (iii) of Lemma 3.3 are
satisfied with $p/q=\hat\tau=d/(R-1)$.

\noindent Then Lemma 3.3 can be
rephrased as
\proclaim Proposition 3.4.  Let $\hat\tau$ be as in Lemma 3.3.
Then
$$\Ztau=\bigcup_{\tuple\in I_\tau}
{\cal Z}\tuple\;.\leqno(3.6)$$

\noindent {\bf Example 3.5.}   In the case of rank two the only
possibility for split bundles is $E_\phi\oplus E_s$, where
$E_s$ is a line bundle of degree $\hat\tau$ and $\phi\in
H^0(E_\phi)$, where $E_\phi$ is a line bundle of degree
$d-\hat\tau$.  The pair $(E_\phi, \phi)$ is determined up to
equivalence by the divisor class of $\phi$ (see [B2]); hence, the
space of equivalence classes of pairs $(E_\phi, \phi)$ is simply the
$d-\hat\tau$ symmetric product $\hbox{Sym}^{d-\hat\tau}\Sigma$.
Since $E_s$ is arbitrary, we have
$${\cal Z}_\tau\simeq \hbox{Sym}^{d-\hat\tau}\Sigma\times {\cal
J}_{\hat\tau}\;,\leqno(3.7)$$
where ${\cal J}_{\hat\tau}$ denotes the component of the Jacobian
variety of $\Sigma$ corresponding to line bundles of degree
$\hat\tau$.

\beginsection 4. The algebraic stratification of $\Bhat$

In this and the next section we
examine two stratifications of the master space
$\Bhat$.  From
one point of view these are a consequence of the fact that the
holomorphic pairs in $\Bhat$
can be characterized by a stability property, and
as such admits two natural
filtrations.
The filtrations are the
analogs of the Seshadri filtration for semistable bundles,
and just as in that case, can be used to
associate ``gradings" to the stable
pairs.  These lead to
stratifications according to grading type.  From a
different perspective,
the gradings and filtrations can be understood in terms
of the Morse
theory of the moment map $f$ on $\Bhat$.
In that context the gradings correspond to the
critical points at infinity on flow lines
either up or down the gradient of the
moment map.  The stratifications thus
correspond to the stratifications given
by the stable (or unstable) manifolds in
the sense of Morse Theory.  This will
be explained in the next section.
We first give a purely
algebraic description.  Let $(E,\phi)$\ be a holomorphic
pair, where $E$\ is here understood as a holomorphic bundle, i.e. the
underlying smooth bundle with a $\dbare$-operator. The two
filtrations are characterized by the parameters
$\mu_+(E)$\ and $\mu_-(E)$,  where
$$\eqalign{\mu_+(E)&=\hbox{max }
\bigl\{\mu(E') : E'\subset E\hbox{ is a holomorphic
subbundle }\bigr\}\cr
\mu_-(E,\phi)&=\hbox{min }\bigl\{\mu(E/E^{\prime\prime}) : E^{\prime
\prime}\subset E
\hbox{ is a  holomorphic
subbundle {\it and} } \phi\in H^0(E^{\prime\prime})\bigr\}.\cr}$$
Generalizing Definition 1.1, we call the pair $(E,\phi)$
stable if
$$\mu_+(E)<\mu_-(E,\phi).$$
This is clearly equivalent to the pair being $\tau$-stable for all
$\mu_+(E)<\tau < \mu_-(E,\phi)$.  Thus the set of isomorphism classes of stable
pairs is the union over all $\tau$\ of the $\Btau$, or equivalently, the union
of the reduced level sets $f^{-1}(\hat \tau)/\U$\ (cf. Proposition 2.18).
Conversely, all  points in $\Bhat$\ are represented by stable pairs.   We will
refer to the two filtrations associated to a stable pair as the
$\mu_-$-filtration and the  $\mu_+$-filtration.

\proclaim Proposition 4.1. {\tenrm (The $\mu_-$-filtration)}
Let $(E,\phi)$\ be a stable holomorphic pair.  There is a
filtration of $E$ by  subbundles
$$0\subset E_\phi=F_0\subset F_1\subset F_2\subset\cdots\subset F_n=E
\leqno(4.1)$$
such that the following properties hold:
\itemitem{(i)} $\phi\in H^0(E_\phi)$,
the pair $(E_\phi,\phi)$ is a stable pair, and
$\mu_+(E_\phi)<\mu_-(E,\phi)<\mu_-(E_\phi,\phi)\;,$
\itemitem{(ii)} for $i=1,\dots,n$
the quotients $F_i/F_{i-1}$ are stable bundles each
of slope $\mu(F_i/F_{i-1})$\break
$=\mu_-(E,\phi)$,
\itemitem{(iii)} $E_\phi$ has minimal rank
among filtrations satisfying (i) and (ii).

\noindent {\sl The
subbundle $E_\phi$ is uniquely determined, and the graded object
$$gr^-(E,\phi) = E_\phi\oplus F_1/F_0\oplus F_2/F_1\oplus
\cdots\oplus E/F_{n-1}$$
is unique up to isomorphism of $F_1/F_0\oplus F_2/F_1\oplus
\cdots\oplus  E/F_{n-1}$. }

Using this result, we define
the $\mu_-$-grading for the pair $(E, \phi)$ by

\proclaim Definition 4.2. The $\mu_-$-grading
for a stable pair $(E, \phi)$  is given by
$$gr^-(E,\phi) = (gr^-(E),\phi)$$
where $gr^-(E)$\ is as above.

We obtain a convenient way
to interpret this grading if we adopt the convention
that the slope of a stable \it pair \rm is
$\mu(E, \phi)=\tau$\ where $\tau$\ is any number such that
$\mu_+(E)<\tau<\mu_-(E,\phi)$.  Then the
pair $(E_\phi,\phi)$\ has slope $\mu_-(E,\phi)$.  By rewriting
$$gr^-(E,\phi)=(E_\phi,\phi)\oplus F_1/F_0\oplus F_2/F_1
\oplus\cdots\oplus  E/F_{n-1},$$
this then becomes a direct sum of stable objects all of the same slope.

In a similar way, the second
filtration will be characterized by the fact that
all quotients are stable objects of
slope $\mu_+(E)$. Such a filtration is well
known in the case when $\mu_+(E)=\mu(E)$, i.e. when $E$
is a semistable  bundle.  Indeed
the usual Seshadri filtration has this property.
If  the bundle is not semistable, and
thus $\mu_+(E)>\mu(E)$, we obtain the  filtration a follows:

\proclaim Proposition 4.3. {\tenrm (The $\mu_+$-filtration)}
Let $(E,\phi)$ be a stable holomorphic pair.
There is a filtration of $E$ by  subbundles
$$0= F_0\subset F_1\subset F_2\subset\cdots\subset
F_n\subset F_{n+1}=E\leqno(4.2)$$
such that the following properties hold:
If $E$ is semistable then this is a
Seshadri filtration and the
quotients $F_i/F_{i-1}$\ are all stable
bundles of slope $\mu_+(E)=\mu(E)$.  Otherwise,
\itemitem{(i)} for $i=1,\ldots, n$ the quotients
$F_i/F_{i-1}$ are stable bundles each
of slope $\mu(F_i/F_{i-1})$\break
$=\mu_+(E)$,
\itemitem{(ii)} $\phi$\ has a non-zero projection,
$\varphi$, into $H^0(E/F_n)$, and
the pair $(E/F_n,\varphi)$ is a stable pair with
$\mu_+(E/F_n)<\mu_+(E)<\mu_-(E/F_n)$,
\itemitem{(iii)} $E/F_n$ has minimal rank among
filtrations satisfying (i) and (ii).

\noindent {\sl In the case
where $\mu_+(E)>\mu(E)$, the quotient $Q=E/F_n$ is uniquely
determined, and the graded object
$$gr^+(E)=F_1/F_0\oplus F_2/F_1\oplus\dots F_n/F_{n-1}\oplus Q$$
is unique up to isomorphism of $F_1/F_0\oplus F_2/F_1\oplus\dots
F_n/F_{n-1}$.}

\proclaim Definition 4.4. For a stable pair $(E,\phi)$ for which
$\mu_+(E)>\mu(E)$, the $\mu_+(E)$-grading is defined to be
$$gr^+(E, \phi)=(gr^+(E),\varphi).$$
For $\mu_+(E)=\mu(E)$, we set
$$gr^+(E,\phi)=(Gr(E),0),$$
where $Gr(E)$ is the grading for $E$ coming from the
Seshadri filtration.

Notice that the filtration for the semistable
subbundle $F_n$ is precisely
the Seshadri filtration.
The case $\mu_+(E)=\mu(E)$ thus corresponds
to the case where $Q=0$.
We now prove Propositions 4.1 and 4.3.  We begin with the
$\mu_+(E)$-filtration.

\proclaim Proposition 4.5.  Given a stable pair  $(E,\phi)$ there is a
unique quotient $Q$ of $E$ arising from an exact sequence
$$0\longrightarrow F\longrightarrow E\longrightarrow Q
\longrightarrow 0\; ,$$
with the properties:
\itemitem{(i)} $F$ is a semi-stable bundle,
\itemitem{(ii)} $\mu(F)=\mu_+(E)$,
\itemitem{(iii)} if $Q\ne 0$,
then under projection of $E$ onto $Q$ the section $\phi$
has a nontrivial image, $\varphi$, and the holomorphic pair
$(Q, \varphi)$ is  stable.
\itemitem{(iv)} If $Q\neq 0$, then $\mu_+(Q)<\mu_+(E)<\mu_-(Q)$,
\itemitem{(v)} $Q$ has minimal rank among  quotients
satisfying (i) - (iv).

\noindent {\it Proof.}
If $E$ is a semistable bundle then $\mu_+(E)=\mu(E)$, and we take
$F=E$, $Q=0$.  Otherwise, for $F$ we take
the unique maximal semistable
subbundle of $E$.  Properties (i) , (ii) and (v) follow immediately from this
choice of $F$. Properties (iii) and (iv)
are consequences of the following

\proclaim Lemma 4.6. Let $(E,\phi)$
be a stable pair with $\mu_+(E)>\mu(E)$.
Let $F$ be the unique maximal semistable subbundle of
$E$, and let $Q$ be the
quotient $E/F$.  Let $\varphi\in H^0(Q)$
be the image of $\phi$ under the
projection of $E$ onto $Q$.  Then
\itemitem{(i)} $\varphi\ne 0$,
\itemitem{(ii)} $\mu_+(Q)<\mu_+(E)$,
\itemitem{(iii)} $\mu_-(Q,\phi)\ge\mu_-(E,\phi)$.

\noindent {\it Proof.}
(i)  If $\varphi=0$ then $\phi\in H^0(F)$.
But then, as $(E,\phi)$\ is stable,
$\mu(E/F)\geq\mu_-(E,\phi)>\mu_+(E)=\mu(F)$, i.e. $\mu(Q)>\mu(F)$.
This is
incompatible with $\mu(F)=\mu_+(E)>\mu(E)$. Thus $\varphi\ne 0$.
\par\noindent
(ii)  Let $Q'\subset Q$ be any holomorphic subbundle.  Lift $Q'$ to a
subbundle $E'\subset E$.  This gives a short exact sequence
$$0\longrightarrow F\longrightarrow E'\longrightarrow Q'
\longrightarrow 0\ .$$
By definition, $\mu(E')\le\mu_+(E)=\mu(F)$,
but in fact the inequality must be
strict since $\rank(E')>\rank(F)$.
It follows from this and the above short
exact sequence that
$\mu(Q')<\mu_+(E)$.  Thus $\mu_+(Q)<\mu_+(E)$.
\par\noindent
(iii)  Suppose in addition that $\varphi\in Q'$.
Then $\phi\in E'$, and thus
$\mu(E/E')\ge\mu_-(E,\phi)$.
But $\mu(E/E')=\mu(Q/Q')$, and thus it follows that
$\mu_-(Q,\phi)\ge\mu_-(E,\phi)$.

\par\noindent
{\it Proof of Proposition 4.3.}
The Proposition follows immediately from Proposition 4.5,
plus the usual  Seshadri filtration for the subbundle $F$.

We now turn to the proof of Proposition 4.1.  The key result here is

\proclaim Lemma 4.7.  Let $(E,\phi)$ be a stable pair.
Let $E_\phi\subset E$
be a holomorphic subbundle such that $\phi\in H^0(E_\phi)$\ and $\mu(E/E_\phi)
= \mu_-(E,\phi)$.  Then
\itemitem{(i)}  $\mu_+(E_{\phi})\leq\mu_+(E)$,
\itemitem{(ii)} $\mu_-(E_{\phi},\phi)\geq\mu_-(E,\phi)$,
and the inequality is strict if $E_\phi$
has minimal rank among all subbundles satisfying the
hypotheses of the Lemma,
\itemitem{(iii)} $(E_\phi,\phi)$ is a stable pair,
\itemitem{(iv)} $E/E_\phi$ is a semi-stable bundle,
\itemitem{(v)} $\mu(E_{\phi})< \mu_-(E,\phi) $
\itemitem{(vi)} Suppose that $E_\phi$
has minimal rank among all subbundles satisfying
the hypotheses of the Lemma, and that $E_\phi '$ is any other subbundle such
that $\phi\in H^0(E'_\phi)$ and $\mu(E/E'_\phi) = \mu_-$.  Then
$E_\phi\subseteq E'_{\phi}$.

\noindent{\it Proof.}
(i)  The first inequality is clear, since $E_{\phi}$
is a subbundle of $E$.
\par\noindent
(ii) Let $E^{\prime\prime}$ be such that $E^{\prime\prime}
\subset E_{\phi}\subset E$,
and $\phi\in E^{\prime\prime}$. Use
the following notation:
$$\matrix{
& E^{\prime\prime} & E_{\phi} & E\cr
\hbox{degree} & d^{\prime\prime} & d_{\phi} & d\cr
\hbox{rank} & R^{\prime\prime} & R_{\phi} & R\cr
}$$
Then
$$\mu(E_{\phi}/E^{\prime\prime})-\mu(E/E^{\prime\prime})=
{R(d_{\phi}-d)+R^{\prime\prime}(d-d_{\phi})+R_{\phi}(d^{\prime\prime}
-d) \over
(R-R^{\prime\prime})(R_{\phi}-R^{\prime\prime})}\;,$$
and
$$\mu(E/E^{\prime\prime})-\mu(E/E_{\phi})=
{R(d_{\phi}-d)+R^{\prime\prime}(d-d_{\phi})+R_{\phi}(d^{\prime\prime}
-d)\over(
R-R^{\prime\prime})(R-R_{\phi})}\;.$$
Hence,
$$\mu(E_{\phi}/E^{\prime\prime})-\mu(E/E^{\prime\prime})=
\left({R-R_{\phi}\over R_{\phi}-R^{\prime\prime}
}\right)\left(\mu(E/E^{\prime\prime})-\mu(E/E_{\phi})\right)\ .$$
The right hand side of this equality is non-negative
by the definition of
$\mu_-(E,\phi)$, and it is strictly positive if $E_\phi$ has minimal rank among
all
subbundles satisfying the hypotheses of the Lemma. This is because
$\mu(E/E^{\prime\prime})=\mu(E/E_{\phi})$ would imply that $E^{\prime
\prime}$
is a subbundle satisfying
the hypotheses but with rank less than that of $E_\phi$.
The result now follows
from the fact that $\mu(E/E_{\phi})=\mu_-(E,\phi)$.
\par\noindent
(iii)  This follows immediately from (i) and (ii), and the fact that
$(E,\phi)$ is  stable.
\par\noindent
(iv)  Suppose that $E/E_\phi$ is not semistable.
Pick a subbundle $F\subset
E/E_\phi$ such that
$\mu(F)=\mu_+(E/E_\phi)$.  Let $E'\subset E$\ be the lift
of $F$ to $E$, i.e. such that
$$0\longrightarrow E_\phi\longrightarrow E'\longrightarrow F
\longrightarrow 0\ .$$
Now $\mu(E/E')=\mu((E/E_\phi)/F)$, and if $\mu(F)>\mu(E/E_\phi)$ then
$\mu((E/E_\phi)/F)<\mu(E/E_\phi)$ .  Hence
$$\mu(E/E')<\mu(E/E_\phi) .$$
However, since $\phi\in H^0(E')$, we have
$$\eqalign{
\mu(E/E')&\ge\mu_-(E,\phi)\cr
&=\mu(E/E_\phi)\;.\cr}$$
Thus $E/E_\phi$ must be semistable.
\par\noindent
(v)  Since $(E,\phi)$ is stable, we have
$\mu(E_\phi)\leq\mu_+(E)<\mu_-(E,\phi)$.
\par\noindent
(vi)  Let $E_\phi$ and $E_\phi '$
satisfy the hypotheses, and suppose that
$E_\phi$ is of minimal rank. Now consider the map
$E_{\phi}\longrightarrow
E/E_{\phi}'$, and let $K$ and $L$ be its kernel and image
respectively. We thus have
$$0\longrightarrow K\longrightarrow E_{\phi}\longrightarrow
L\longrightarrow 0\leqno(4.3)$$
Suppose that $K\neq E_\phi$.
Since $\phi$ is a section of $K$, we have
$$\mu(E_{\phi}/K)\geq \mu_-(E_{\phi},\phi)\;. $$
Also, since by (iv) $E/E_{\phi}'$ is semistable, we have
$$\mu(L)\leq\mu(E/E_{\phi}')=\mu_-(E,\phi) \;.$$
By (ii) we have $\mu_-(E,\phi)<\mu_-(E_{\phi},\phi)$,
and from (4.3) we have
$\mu(E_{\phi}/K)=\mu(L)$.  We thus get
$$\mu(L)\leq\mu_-(E,\phi)<\mu_-(E_{\phi},\phi)
\leq\mu(E_{\phi}/K)=\mu(L)\;,$$
which is impossible.  We conclude that
$K=E_\phi$, i.e. $E_\phi\subset
E'_{\phi}$.

\proclaim Proposition 4.8. Given a stable pair $(E,\phi)$
there is a unique
subbundle $E_\phi\subset E$ such that
\itemitem{(i)}
$\phi\in H^0(E_\phi)$\ and $(E_\phi,\phi)$ is a stable pair,
\itemitem{(ii)} $E/E_\phi$ is a semi-stable bundle,
\itemitem{(iii)} $\mu(E/E_\phi)=\mu_-(E,\phi)$,
\itemitem{(iv)} $\mu_+(E_{\phi})< \mu_-(E,\phi)<\mu_-(E_{\phi},\phi) $
\itemitem{(v)} $E_\phi$ has minimal rank among all subbundles satisfying
(i)-(iv)

\noindent {\it Proof.}
By parts (i)-(v) of Lemma 4.7, any subbundle $E_\phi\subset E$
such that
$\phi\in H^0(E_\phi)$ and $\mu(E/E_\phi) = \mu_-$ will satisfy
(i)-(iv).  By
part (vi) of Lemma 4.7, there is a unique such $E_\phi$ of minimal rank.

\noindent {\it Proof of Proposition 4.1.}
The required filtration is constructed as
follows.  Let
$$0\subset Q_1\subset Q_2\subset\cdots\subset Q_n=E/E_\phi$$
be the Seshadri filtration for $E/E_\phi$.
Set $F_i=\pi^{-1}(Q_i)$,
where $\pi:E\longmapsto E/E_\phi$ is the projection map.

\noindent {\it Remark.} An
important feature of the two gradings, $gr^-(E,\phi)$ and
$gr^-(E,\phi)$
is that, as pairs, they are both semistable and thus correspond
to points in the masterspace $\Bhat$.  Indeed for $gr^-(E,\phi)$
we have
$\mu_-(gr^-(E))=\mu_+(gr^-(E))=\mu_-(E)$,
while for $gr^+(E,\phi)$ we have
$\mu_-(gr^+(E))=\mu_+(gr^+(E))=\mu_+(E)$.

It is worth pointing out that there
are other filtrations  associated to a stable pair,
but none of the resulting
gradings are semistable as pairs and thus are not represented in
$\Bhat$.
These other gradings arise by successive
applications of Propositions 4.5 and
4.8.  For example, one can successively apply
Proposition 4.8 to $E_\phi$ in
the extension
$$0\longrightarrow E_\phi \longrightarrow  E\longrightarrow Q
\longrightarrow 0\;,$$
until the subbundle containing $\phi$
is of rank one.  This yields a grading of  the form
$([\phi]\oplus Gr(E/[\phi]), \phi)$, where $[\phi]$
is the line subbundle
generated by $\phi$, and $Gr(E/[\phi])$
is the Seshadri grading for
$E/[\phi]$.  Similarly, successive
application of Proposition 4.5 to the
quotient pair $(Q,\varphi)$
leads in all cases (i.e. not just when $E$ is
semistable) to the grading $(Gr(E),0)$.
Such gradings and their relation to
$gr^{\pm}(E, \phi)$ will be discussed in a future publication.
We now use the gradings $gr^\pm(E,\phi)$ to define a stratification of
$\Bhat$.
Notice firstly that for any pair $(\dbf)$, the gradings
$gr^\pm(E,\phi)$
are each characterised by (n+2)-tuples of integers $\tuple$
which satisfy the constraints in Lemma 3.3.  The notation is
such that the pair
$(d_{\phi},R_{\phi})$ refers to the degree and rank of the summand in
$gr^\pm(E,\phi)$\ which contains the section.
In the case of $gr^-(E,\phi)$
this is $E_{\phi}=F_0$ in the notation of  (4.1),
while for $gr^+(E,\phi)$
this is the quotient $Q$ of Proposition 4.5. In both cases the $R_i$
give the
ranks of the quotients $F_i/F_{i-1}$.

\proclaim Definition 4.9.
Given a stable pair $(E,\phi)$ and an $(n+2)$-tuple
$\tuple\in I_\tau$, let
$$\eqalign{{\cal W}^\pm\tuple &=
\Bigl\{ (E,\phi)\in\Bhat : (E,\phi) \hbox{ is a stable pair,
and}\cr
gr^\pm(E,\phi)&\in{\cal Z}\tuple
\Bigr\}\,\bigcup\, {\cal Z}\tuple\;.\cr}$$
Set
$$\Wtau^\pm=
\bigcup_{\tuple\in I_\tau}
{\cal W}^\pm\tuple\;.$$

The next proposition justifies our definition of ${\cal W}^\pm\tuple$.
Its proof is straightforward, and since we shall not have need of
the statement in the sequel we omit the details.

\proclaim Proposition 4.10.  With respect to the obvious ordering of
critical values of $\hat\tau$ in \break
$[d/R, d/(R-1)]$, the subspaces $\Wtau^+$
form a piecification of $\Bhat$ in the sense of Goresky-MacPherson
(see [G-M]).  A similar result holds for $\Wtau^-$.  Moreover, the
subspaces
$${\cal W}^\pm\tuple \cap\Bhat_0$$
are V-manifolds for all $\tuple\in I_\tau$.

We end this section with the following proposition which will be
used in \S 6:

\proclaim Proposition 4.11.  If $R>2$, then for critical values
$\hat\tau \in (d/R, d/(R-1))$ the complex codimension of
$\Wtau^\pm$ in $\Bhat$ is $\geq 2$.

\noindent For the proof we shall need the following simple
\proclaim Lemma 4.12.  Suppose $(E_\phi, \phi)$ is a $\tau$-stable
pair and $E_s$ is a semistable bundle with slope $\tau$.  Then
$H^0(E_\phi\otimes E_s^\ast)=\{0\}$.

\noindent {\it Proof.}  Suppose $\alpha\in
H^0(E_\phi\otimes E_s^\ast)$ and $\alpha\not\equiv 0$.  Then
$\alpha$ defines a map of sheaves $E_s\to E_\phi$, and since
$\alpha\not\equiv 0$, $\rank(\ker\alpha) < \rank(E_s)$.  But then
the semistability of $E_s$ and the $\tau$-stability of $E_\phi$
imply
$$\tau=\mu(E_s)\leq \mu(E_s/\ker\alpha)=\mu(\hbox{image }\alpha)
<\tau\;,$$
which  is a contradiction.  This proves the Lemma.

\noindent {\it Proof of Proposition 4.11. }
It suffices to compute the codimension of
the largest stratum ${\cal W}^\pm(d_\phi, R_\phi, R_s)$, where
$(d_\phi, R_\phi, R_s)\in I_\tau$.  We first consider ${\cal W}^-$.
Let $(E, \phi)$ be a stable pair such that
$gr^+(E,\phi)\in {\cal Z}(d_\phi, R_\phi, R_s)$.
Then we have an exact  sequence
$$0\mapright{} E_\phi\mapright{} E\mapright{} E_s\mapright{} 0\;.$$
The tangent space to ${\cal W}^-(d_\phi, R_\phi, R_s)$ at $(E,\phi)$
naturally splits
$$T_{(E,\phi)}{\cal W}^-(d_\phi, R_\phi, R_s)\simeq T_{(E_\phi,\phi;
E_s)}{\cal Z}(d_\phi, R_\phi, R_s)\oplus \hbox{Ext}^1(E_s,
E_\phi)\;.$$
The dimension of ${\cal Z}(d_\phi, R_\phi, R_s)$ is computed as in
Section 3 of [B-D1]:
$$\dim{\cal Z}(d_\phi, R_\phi, R_s)=d-\hat\tau R_s
+1+(R-R_s-1)(R-R_s)(g-1)+ R_s^2(g-1)\;,$$
and $\hbox{Ext}^1(E_s, E_\phi)\simeq H^1(E_\phi\otimes E_s^\ast)$.
By Lemma 4.12 and Riemann-Roch we have that

$$\eqalign{\dim T_{(E,\phi)}{\cal W}^-(d_\phi, R_\phi, R_s)
&= \dim{\cal Z}(d_\phi, R_\phi, R_s)
+\left( \hat\tau R -d + (R-R_s)(g-1)\right) R_s\cr
&= d+1+
(R^2-R)(g-1)+R_s(R_s-R+1)(g-1)\cr
&\qquad +\left(\hat\tau(R-1)-d\right)
R_s\;.\cr}$$
By Theorem 2.16,
$$p_-(R_s,\hat\tau)\equiv
\hbox{codim } {\cal W}^-(d_\phi, R_\phi, R_s)
=\left(d-\hat\tau(R-1)\right)R_s+ R_s(R-R_s-1)(g-1)\;.$$
For $1< R_s < R-1$ the last term in the expression above is $\geq 2$
(we assume $g>1$), and since $\hat\tau < d/(R-1)$ the first term is
positive.  Now we check the case where $R_s=1$. Then
$$p_-(1,\hat\tau)= d-\hat\tau(R-1)+(R-2)(g-1)\;.$$
Since we assume $R>2$ and $g>1$ the last term is $\geq 1$.  Also,
$d-\hat\tau(R-1)$ is a positive integer and so must also be $\geq
1$.  Thus, $p_-(1,\hat\tau)\geq 2$.  For $R_s=R-1$,
$$p_-(R-1, \hat\tau)=\left( d-\hat\tau(R-1)\right)(R-1)\;.
\leqno(4.4)$$
Again, $d-\hat\tau(R-1)$ is a positive integer and $R-1\geq 2$.
Therefore, in all cases we have
$\hbox{codim }{\cal W}^-(d_\phi, R_\phi, R_s)\geq 2$.
Now consider
${\cal W}^+(d_\phi, R_\phi, R_s)$.
Let $(E, \phi)$ be a stable pair such that $gr^+(E,\phi)\in
{\cal Z}(d_\phi, R_\phi, R_s)$.
Then $E$ may be written (see Proposition
4.3):
$$0\mapright{} F\mapright{} E\mapright{\pi} Q\mapright{} 0\;.$$
As in the case of  ${\cal W}^+$  the tangent space
$$T_{(E,\phi)}{\cal W}^+(d_\phi, R_\phi, R_s)\supset
T_{(Q, \pi(\phi); F)}{\cal Z}(d_\phi, R_\phi, R_s)$$
as a summand.  The complement is naturally isomorphic to the space
of extensions of $Q$ by $F$ direct sum with the equivalence
classes of liftings of $\pi(\phi)$.  The liftings are parameterized
by $H^0(F)$, and two liftings are equivalent if and only if
they differ by an element of $H^0(Q^\ast\otimes F)$.
Therefore,
$$\eqalign{
\dim T_{(E,\phi)}{\cal W}^+(d_\phi, R_\phi, R_s)&=
\dim T_{(Q, \pi(\phi); F)}{\cal Z}(d_\phi, R_\phi, R_s)
+\dim H^1(Q^\ast\otimes F)\cr
&\qquad +\dim H^0(F)-\dim
H^0(Q^\ast\otimes F)\;.\cr}$$
Since $F$ is stable with slope $\hat\tau > d/R > 2g-2$, $H^1(F)=0$.
Therefore, by Riemann-Roch
$$\eqalign{
\dim T_{(E,\phi)}{\cal W}^+(d_\phi, R_\phi, R_s)&=
d-\hat\tau R_s +1 +(R-R_s-1)(R-R_s)(g-1)+ R_s^2(g-1)\cr
&\qquad +\left(d-\hat\tau(R-1)+(R-R_s-1)(g-1)\right) R_s\cr
&= d+1+ (R^2-R)(g-1) +(d-R\hat \tau)R_s\cr
&\qquad +R_s(R_s-R)(g-1)\;.\cr}$$
By Theorem 2.16, this implies
$$p_+(R_s,\hat\tau)\equiv
\hbox{codim } {\cal W}^+(d_\phi, R_\phi, R_s)
=(R\hat\tau-d)R_s +R_s(R-R_s)(g-1)\;.\leqno(4.5)$$
Since $\hat\tau > d/R$,
$$p_+(R_s,\hat\tau) > R_s(R-R_s)(g-1)\;,$$
and it is easily checked that the latter expression is always $\geq
2$ for $R>2$, $g>1$, and $1\leq R_s\leq R-1$.  This completes the
proof of Proposition 4.11.

\beginsection 5. The Morse theory of $f$

We now turn to the description of the Morse theory of the function
$f$.  We shall write down solutions to the gradient flow of $f$ and
describe the stable and unstable manifold stratifications of
$\Bhat$.  Furthermore, we show that the Morse theoretical
stratification of $\Bhat$ coincides with the algebraic
stratification of the previous section.  The results of this section
are similar in spirit to the results of [D].  However, the situation
here is technically simpler because we are dealing with a finite
dimensional problem and an abelian group action (see also [K1]).

\proclaim Proposition 5.1.  Let $\Phi:\Bhat\times
[0,\infty)\to\Bhat$ be the flow
$$\Phi_t[\dbf]=[\dbare , e^{-t/2\pi R}\phi]\;.$$
Then $\Phi$ is continuous.  Moreover, $\Phi$ preserves $\Bhat_0$ and
coincides with the gradient flow of $f$ on $\Bhat_0$.

\noindent {\it Proof.}  We must verify that
$${d\Phi_t\over dt}=-\nabla_{\Phi_t} f\;.\leqno(5.1)$$
First recall that after identifying $T_{[\dbf]}\Bhat$ with
$H^1(\Cxo)$, the infinitesimal vector field of the $\U$ action on
$\Bhat$ is given by
$$\xi^\#[\dbf]={i\over R}(0,\phi)\;.\leqno(5.2)$$
Indeed,

$$\xi^\#[\dbf]=\left.{d\over d\theta}\right|_{\theta=0}
[\dbare, e^{i\theta/R}\phi]=
{i\over R}(0,\phi)\;.$$
Moreover,
$$\eqalign{\nabla_{\Phi_t[\dbf]} f &=
{-1\over 2\pi i}\nabla_{\Phi_t[\dbf]} \Psi
={1\over 2\pi i}\xi^\# (\Phi_t[\dbf])\cr
&= {1\over 2\pi R}(0, e^{-t/2\pi R}
\phi)=-{d\over dt}(0, e^{-t/2\pi R}\phi)
=-{d\Phi_t[\dbf]\over dt}\;,\cr}$$
which is what was to be shown.

\proclaim Definition 5.2. Given a critical $\hat\tau$ and

$\tuple\in I_\tau$ as in \S 3, let
$${\cal W}^s\tuple =\left\{ [\dbf]\in\Bhat :
\lim_{t\to\infty}\Phi_t[\dbf]\in {\cal Z}\tuple\right\}\;,$$
and let ${\cal W}^u\tuple$ be defined similarly as $t\to -\infty$.
Also, for a critical value of $\hat\tau$, we set
$$\eqalign{
\Wtau^s &=\bigcup_{\tuple\in I_\tau}\overline{\cal W}^s\tuple\cr
\Wtau^u &=\bigcup_{\tuple\in I_\tau}\overline{\cal W}^u\tuple\;.\cr}$$
We call $\{\Wtau^s\}$ and $\{\Wtau^u\}$ the stable and unstable
Morse stratifications of $\Bhat$, respectively.

\proclaim Theorem 5.3.  For each critical value $\tau$,
$\Wtau^s=\Wtau^+$, and $\Wtau^u=\Wtau^-$.  Consequently, the Morse
stratification of $\Bhat$ coincides with the algebraic stratification
of \S 3.

\noindent {\it Proof.} We shall show that $\Wtau^u=\Wtau^-$.
Indeed, since both $\{\Wtau^u\}$ and $\{\Wtau^-\}$ are
stratifications of $\Bhat$, it suffices to prove the inclusion
$\Wtau^-\subset \Wtau^u$ for all $\tau$.  In fact, we are going to
show that
$${\cal W}^-\tuple\subset {\cal W}^u\tuple$$
for all $\tuple\in I_\tau$.  Fix $[\dbf]\in{\cal W}^-\tuple$.
Let
$$0=E_\phi = F_0\subset F_1\subset\cdots\subset F_n=E$$
denote the $\mu_-$ filtration of the pair $(E,\phi)$.  Fix real
numbers
$$0 < \mu_1 < \mu_2 < \cdots < \mu_n \;,
\qquad \sum_{i=1}^n R_i\mu_i =R_\phi\;,$$
and consider the following
1-parameter subgroup of gauge transformations in $\Gco$,
$$ g_t=\left(
\matrix{e^{t/2\pi R}&0&\cdots&0\cr
0&e^{-t\mu_1/2\pi R}&\cdots&0\cr
\vdots&\vdots&\ddots&\vdots\cr
0&0&\cdots&e^{-t\mu_n/2\pi R}\cr}\right)\leqno(5.3)$$
written diagonally with respect to the filtration above.
Then
$$\eqalign{\lim_{t\to -\infty}\Phi_t[\dbf] &=
\lim_{t\to -\infty}[\dbare, e^{-t/2\pi R}\phi]\cr
&=\lim_{t\to -\infty}[\dbare, g_t^{-1} \phi]\cr
&=\lim_{t\to -\infty}[g_t(\dbare), \phi]\cr
&=\left[gr^-(E,\phi), \phi\right]\;.\cr}$$
The last equality follows the same way as in [D], p. 716.
Hence,
$$[\dbf]\in{\cal W}^u\tuple\;.$$
The case of the stable
manifolds is similar.
To prove that ${\cal W}^+\tuple\subset{\cal W}^s\tuple$ for all
$\tuple\in I_\tau$ one must show that for all \break
$[\dbf]\in{\cal
W}^+\tuple$,
$$\lim_{t\to\infty}\Phi_t[\dbf]=[gr^+(E,\phi),\varphi]\;,$$
where $[gr^+(E,\phi), \varphi]$ is the $\mu_+$-grading as defined in
\S 4.  The above method can be used, but now the complex gauge
transformations $g_t\in\Gco$ must be defined as follows.
Let
$$0= F_0\subset F_1\subset\cdots\subset F_{n+1}=E$$
denote the $\mu_+$-filtration of the pair $(E,\phi)$.  If $E$ is a
semistable bundle, fix real numbers
$$1>\mu_1 > \mu_2 > \cdots > \mu_{n+1}\;,
\qquad \sum_{i=1}^n R_i\mu_i =0\;,$$
and let $g_t$ be
$$ g_t=\left(
\matrix{e^{t\mu_1 /2\pi R}&0&\cdots&0\cr
0&e^{t\mu_2/2\pi R}&\cdots&0\cr
\vdots&\vdots&\ddots&\vdots\cr
0&0&\cdots&e^{t\mu_{n+1}/2\pi R}\cr}\right)\leqno(5.3)$$
written diagonally with respect to the filtration above.
If $E$ is not semistable, then one must take $\mu_{n+1}=1$ and impose
the constraint $R_{n+1}
+ \sum_{i=1}^n R_i\mu_i =0$.
The rest of the argument proceeds as before.

\beginsection 6.  Birational equivalence of stable pairs

In this section we describe how the moduli of vortices $\Btau$
change with respect to $\tau$.  The analogous situation has been
studied in the symplectic category by Guillemin and Sternberg [G-S]
and in the algebraic category by Goresky and MacPherson [G-M].
 However, since  $\Bhat$ has
singularities and  no obvious embedding in projective space
compatible with the $\U$ action, the results of [G-S] and [G-M] are
not directly applicable to the case at hand.  We thus prove the
following directly:

\proclaim Theorem 6.1.  (i)  Suppose the interval
$[\hat\tau,\hat\tau+\varepsilon]$   contains no critical value of
the function $f$.  Then the Morse flow induces a biholomorphism
between ${\cal B}_{\tau+\varepsilon}$ and $\Btau$.  (ii)  Suppose
that $\hat\tau$ is the only critical value of $f$ in the interval
$[\hat\tau,\hat\tau+\varepsilon]$.  Then the Morse flow defines
a continuous map from ${\cal B}_{\tau+\varepsilon}$ onto $\Btau$ which
restricts to a biholomorphism between ${\cal
B}_{\tau+\varepsilon}\setminus
{\Bbb P}(\Wtau^+)$ and $\Btau\setminus\Ztau$, where
$${\Bbb P}_\varepsilon
(\Wtau^+)=\Wtau^+\cap f^{-1}(\tau+\varepsilon)/\U\;.$$
(iii)  In the case $R=2$, the restriction of the Morse flow to ${\Bbb
P}_\varepsilon(\Wtau^+)$ induces a map
$${\Bbb P}_\varepsilon(\Wtau^+)\mapright{\pi}{\cal Z}_\tau$$
which is a holomorphic projective bundle (unless $d$ is even and
$\hat\tau=d/2$).
In particular, in rank two ${\Bbb P}_\varepsilon(\Wtau^+)$ is a
smooth subvariety of ${\cal B}_{\tau+\varepsilon}$.

\noindent {\it Proof.}  (i) For the sake of notational simplicity we
denote the equivalence class of the pair $[\dbf]\in\Bhat$ by $x$.
Let
$$F: f^{-1}(\hat\tau+\varepsilon)\times[0,\infty)\longrightarrow{\Bbb
R}\leqno(6.1)$$
denote the map $F(x,t)=f(\Phi_t(x))$. By our assumption on
$[\hat\tau,\hat\tau+\varepsilon]$, $F$ is smooth.  Moreover, given
$(x,t)\in F^{-1}(\hat\tau)$ we have
$$\left.{\partial F\over \partial t}\right|_{(x,t)}=
df_{\Phi_t(x)}\left(\left.{\partial\Phi_t\over \partial
x}\right|_x\right) =\Vert \nabla_{\Phi_t(x)} f\Vert^2 \neq 0\;,$$
since $\hat\tau=f(\Phi_t(x))$ is not a critical value of $f$.  Then
by the implicit function theorem we can solve $F(x,t)=\hat\tau$ as
$t=t(x)$, where $t$ is a smooth function of $x$.  We define
$$\hat\sigma_+ : f^{-1}(\hat\tau+\varepsilon)\longrightarrow
f^{-1}(\hat\tau)\leqno(6.2)$$
by $\hat\sigma_+(x)=f(x, t(x))$.  It follows that $\hat\sigma_+$ is
a diffeomorphism between $f^{-1}(\hat\tau+\varepsilon)$ and
$f^{-1}(\hat\tau)$.

Next we show that $\hat\sigma_+$ is a CR-map with respect to the
induced CR-structure on the level sets
$f^{-1}(\hat\tau+\varepsilon)$ and $f^{-1}(\hat\tau)$.  Indeed, let
$$X\in T^{1,0}\Bhat \cap
Tf^{-1}(\hat\tau+\varepsilon)\otimes{\Bbb C}\;,$$
and let $\overline X$ denote the complex conjugate.  Then
$$d\hat\sigma_+(\overline X)={\partial\Phi_t\over \partial
t}\left({\partial t\over \partial x}(\overline X)\right)+{\partial
\Phi_t\over\partial x}(\overline X)=\nabla f\left({\partial
t\over\partial x}(\overline X)\right)+{\partial \Phi_t\over\partial
x}(\overline X)\;.$$
Since $\Phi_t$ is holomorphic in $x$, $\partial\Phi_t/\partial
x(\overline X)=0$.  Hence,

$$d\hat\sigma_+(\overline X)= \nabla f\left({\partial
t\over\partial x}(\overline X)\right) $$
is both tangential and normal to $f^{-1}(\hat\tau)$, and therefore
$d\hat\sigma_+(\overline X)=0$.  Thus, $\hat\sigma_+$ is a CR-map.
Since $\hat\sigma_+$ and the CR-structures on
$f^{-1}(\hat\tau+\varepsilon)$ and $f^{-1}(\hat\tau)$ are
\U-invariant, $\hat\sigma_+$ induces a biholomorphism
$$\sigma_+: {\cal
B}_{\tau+\varepsilon}=f^{-1}(\hat\tau+\varepsilon)/\U\longrightarrow
f^{-1}(\hat\tau)/\U =\Btau\;.\leqno(6.3)$$
(ii)  The same argument as in (i) gives a smooth map
$$\hat\sigma_+ : f^{-1}(\hat\tau+\varepsilon)\setminus
\Wtau^+\longrightarrow
f^{-1}(\hat\tau)\setminus\Ztau\;.\leqno(6.4)$$
We extend $\hat\sigma_+$ across $\Wtau^+$ by setting
$\hat\sigma_+(x)=\lim_{t\to\infty}\Phi_t(x)$ for $x\in\Wtau^+$.  We
are going to show that $\hat\sigma_+$ is continuous.  It is easily
seen (e.g. from the uniqueness of the filtration) that the
restrictions of $\hat\sigma_+$ to $f^{-1}(\hat\tau_\varepsilon)\setminus
\Wtau^+$ and $\Wtau^+$ are continuous.  Therefore, it suffices to
prove that if $\{x_l\}$ is a sequence in
$f^{-1}(\hat\tau+\varepsilon)$, $x\in\Wtau^+$, and $x_l\to x$, then
$\hat\sigma_+(x_l)\to \hat\sigma_+(x)$.  Let $\rho$ be a metric
compatible with the topology of $\Bhat$.  Set $t_l=t(x_l)$, where
$\hat\sigma_+(x_l)=f(\Phi_{t(x_l)}(x_l))$.  Then
$$\rho\left(\hat\sigma_+(x_l),\hat\sigma_+(x)\right)\leq
\rho\left(\Phi_{t_l}(x_l),\Phi_{t_l}(x)\right)+\rho\left(\Phi_{t_l}
(x), \hat\sigma_+(x)\right)\;.$$
Since clearly $t_l\to\infty$ and $\Phi_t$ is unformly continuous,
both terms on the right hand side of the above inequality go to
zero, and this proves the continuity of $\hat\sigma_+$.  Since
$\hat\sigma_+$ is also \U-invariant, it induces a continuous map
$$\sigma_+: {\cal
B}_{\tau+\varepsilon}=f^{-1}(\hat\tau+\varepsilon)/\U\longrightarrow
f^{-1}(\hat\tau)/\U=\Btau\;.\leqno(6.5)$$
On the other hand, by the same argument as in (i), $\sigma_+$
defines a biholomorphism onto its image away from
$${\Bbb P}_\varepsilon({\cal W}^+)={\cal W}^+\cap
f^{-1}(\hat\tau+\varepsilon)/\U\;.$$
(iii) First, we suppose that $\hat\tau > d/2$ since otherwise ${\Bbb
P}_\varepsilon(\Wtau^+)={\cal B}_{\tau+\varepsilon}$.  Then the
fixed point sets ${\cal Z}_\tau$ in $\Bhat$ are smooth, and hence
${\cal W}^+\cap f^{-1}(\hat\tau+\varepsilon)$ is a smooth
submanifold of $\Bhat\setminus f^{-1}(d/2)$.
The Morse flow clearly
induces a continuous map
$${\cal W}^+\cap f^{-1}(\hat\tau+\varepsilon)\mapright{\pi} {\cal
Z}_\tau\;,\leqno(6.6)$$
which is an odd dimensional sphere bundle (say with fiber
$S^{2n+1}$) over ${\cal Z}_\tau$.  Since ${\cal W}^+$ is an analytic
subvariety, the CR-structure on $f^{-1}(\hat\tau+\varepsilon)$
induces a CR-structure on the intersection with ${\cal W}^+$, and as
in the proof of part (ii) above, the map $\pi$ is a CR-map.  Since
$\pi$ is also $\U$-invariant and the $\U$ action is CR, $\pi$
descends to a holomorphic map
$${\cal W}^+\cap f^{-1}(\hat\tau+\varepsilon)/\U \mapright{\pi} {\cal
Z}_\tau\;,\leqno(6.7)$$
with fiber $S^{2n+1}/\U \simeq {\Bbb P}^n$.
This completes the proof of Theorem 6.1.

By reversing the orientation of the flow lines, we obtain
a similar result
relating $\Btau$ and ${\cal B}_{\tau-\varepsilon}$.  Combining the
two results immediately proves the

\proclaim Corollary 6.2.  If $\hat\tau$ is the only critical value
of $\Psi$ is $[\hat\tau-\varepsilon, \hat\tau+\varepsilon]$, then
${\cal B}_{\tau-\epsilon}$, ${\cal B}_{\tau+\varepsilon}$ are
related by the diagram
$$\matrix{{\cal B}_{\tau-\varepsilon}&&{\cal
B}_{\tau+\varepsilon}\cr
{}_{\sigma_-}&\searrow\qquad \swarrow &{}_{\sigma_+}\cr
&\Btau &\cr}$$
where $\sigma_\pm$ are continuous maps.  Moreover,
$$\sigma_\pm : {\cal
B}_{\tau\pm\varepsilon}\setminus\sigma_\pm^{-1}(\Ztau)\longrightarrow
\Btau\setminus\Ztau$$
are biholomorphisms.

Before continuing, we digress to prove that the $\Btau$ are in fact
projective varieties.

\proclaim Theorem 6.3. {\tenrm (see also [B-D2])}
For all non-critical values of $\hat\tau$,
$\Btau$ is a non-singular projective variety.

\noindent {\it Proof.}  According to [M],
since $\Btau$ is a K\"ahler manifold we need only prove that $\Btau$
is Moishezon.  This is equivalent, by a theorem of Siu (see [S],
Theorem 1), to proving that $\Btau$ admits a hermitian, holomorphic
line bundle which is semipositive and positive at at least one
point.  From the codimension estimate given in
(4.4) and (4.5) it suffices to prove this for $\hat\tau $ close to
$d/R$ (for rank two and $d-1 < \hat\tau < d$  we use a separate,
easier argument).  In this case, according to [B-D1], Theorem 6.4,
there is a holomorphic map
$$\Btau\mapright{\pi} {\cal M}(d,R)\;,$$
where ${\cal M}(d,R)$ denotes the Seshadri compactification of
stable bundles.  Moreover, the restriction of $\pi$ to
the open set ${\cal M}^s(d,R)$
consisting of stable bundles is a fibration with fiber ${\Bbb P}^N$:
$$\Btau^s\mapright{\pi} {\cal M}^s(d,R)\;,$$
{}From [B-D2] there is an hermitian, holomorphic line bundle
$\gamma$ on $\Btau^s$ whose restriction to the fiber is ${\cal
O}_{{\Bbb P}^N}(R)$ (an ${\cal O}(1)$ is not always possible, due to
the Brauer obstruction on ${\cal M}(d,R)$).  By pulling back a
sufficiently high power $k$ of an ample bundle $H\to {\cal M}^s(d,R)$
we can arrange $L=\gamma\otimes \pi^\ast H^k$
to be  positive  at a point.  It is easily seen that
 $L$ extends to a semipositive line bundle on $\Btau$
which is positive at a point
of $\Btau^s$.  Siu's theorem then completes the
proof.

We then have the following
\proclaim Corollary 6.4.  For all noncritical values of $\hat\tau$
in $(d/R, d/(R-1))$, the spaces $\Btau$ are all birational.

\noindent {\it Proof.}  According to Corollary 6.2, the complex
manifolds ${\cal B}_{\tau\pm\varepsilon}\setminus\sigma_\pm^{-1}(\Ztau)$
are biholomorphic.  Thus their fields of meromorphic functions
${\cal M}\left({\cal B}_{\tau\pm\varepsilon}\setminus
\sigma_\pm^{-1}(\Ztau)\right)$
are isomorphic.  On the other hand, since for $R> 2$,
$\sigma_\pm^{-1}(\Ztau)$ has codimension at least 2 in ${\cal
B}_{\tau\pm\varepsilon}$ (see Proposition 4.11) and ${\cal
B}_{\tau\pm\varepsilon}$ are smooth, it follows from the
Kontinuit\"atssatz for meromorphic functions (cf. [K-K], 53A.9) that
$${\cal M}\left({\cal B}_{\tau\pm\varepsilon}
\setminus\sigma_\pm^{-1}(\Ztau)\right) \simeq
{\cal M}\left({\cal B}_{\tau\pm\varepsilon}\right)\;.$$
By Theorem 6.3, ${\cal B}_{\tau\pm\varepsilon}$ are
projective varieties, hence by GAGA (cf. [G-H], p. 171)
${\cal M}\left({\cal B}_{\tau\pm\varepsilon}\right)
\simeq {\cal K}\left({\cal B}_{\tau\pm\varepsilon}\right)$, where
${\cal K}\left({\cal B}_{\tau\pm\varepsilon}\right)$ denotes the
field of rational functions.
It follows  that
${\cal K}\left({\cal B}_{\tau-\varepsilon}\right)\simeq
{\cal K}\left({\cal B}_{\tau+\varepsilon}\right)$, and hence ${\cal
B}_{\tau-\varepsilon}$ is birational to ${\cal
B}_{\tau+\varepsilon}$.

It is interesting to apply Corollary 6.4 to the case of ${\cal
B}_{{d\over R}+\varepsilon}$ and ${\cal B}_{{d\over R-1} -\varepsilon}$.
Assume that $d$ is coprime to both $R$ and $R-1$ and $d>R(2g-2)$.  Let
$U(d,R)\to\Sigma\times{\cal M}(d,R)$ denote the universal bundle
over the moduli space of vector bundles of rank $R$ and degree $d$,
and let $\pi:\Sigma\times {\cal M}(d,R)\to{\cal M}(d,R)$ denote the
projection onto the second factor.  It follows from \S 3 that ${\cal
B}_{{d\over R}+\varepsilon}$ is biholomorphic to the
projectivization of the vector bundle $\pi_\ast U(d,R)$.  On the
other hand, let $\hbox{Ext}^1(U(d,R-1), {\cal O})$ denote the
bundle over $\Sigma\times {\cal M}(d, R-1)$ whose fiber over a point
$(p, E^{\dbare})$ consists of the extensions of
$U(d,R-1)\bigr|_{\Sigma\times\{E^{\dbare}\}}$
by $\cal O$. It follows again
from \S 3 that ${\cal B}_{{d\over R-1}-\varepsilon}$ is
biholomorphic to the projectivization of the restriction of
$\hbox{Ext}^1(U(d,R-1),{\cal O})$ to $\{\hbox{point}\}\times{\cal
M}(d, R-1)$.  By combining with Corollary 6.4, we obtain
\proclaim Corollary 6.5.  Assume that $d> R(2g-2)$ is coprime to
both $R$ and $R-1$.  Then ${\Bbb P}\left(\pi_\ast U(d,R)\right)$
over ${\cal M}(d,R)$ is birational to ${\Bbb
P}\left(\hbox{Ext}^1(U(d,R-1),{\cal O})\right)$ over ${\cal
M}(d,R-1)$.

\noindent
Presumably, Corollary 6.5 may also be obtained by carrying out a GIT
construction of these spaces as in [Be] and [Th].

In the case of rank two  our theorem combined with the result
of [G-S] implies the following theorem of Thaddeus [Th]:
\proclaim Theorem 6.6.  Let $R=2$ and $d > 4(g-1)$.  Suppose that
$\hat\tau\in (d/2, d)$ is the only critical value of $f$ in the interval
$[\hat\tau-\varepsilon, \hat\tau+\varepsilon]$.
Then there is a projective variety $\tilde\Btau$ and holomorphic
maps
$$\matrix{&\tilde\Btau &\cr
	{}^{\rho_-}&\swarrow\qquad\searrow &{}^{\rho_+}\cr
	{\cal B}_{\tau-\varepsilon}&&{\cal
	B}_{\tau+\varepsilon}\cr}$$

\noindent {\sl Moreover, for $\hat\tau < d-1$,
$\rho_\pm$ are blow-down maps onto the smooth
subvarieties ${\Bbb P}_\varepsilon(\Wtau^\pm)$.
For $\hat\tau = d-1$, $\rho_+$ is the blow-down map
onto ${\Bbb P}_\varepsilon (\Wtau^+)$ and $\rho_-$ is the identity.
}

\noindent {\it Proof.}
By Theorem 6.1 (iii), ${\Bbb P}_\varepsilon(\Wtau^\pm)$ are smooth
subvarieties of ${\cal B}_{\tau\pm\varepsilon}$, respectively.
Their codimensions are given by
$$\eqalign{p_+(\hat\tau)&= 2\hat\tau -d +g-1\cr
p_-(\hat\tau)&= d-\hat\tau\cr}$$
(see (4.4) and (4.5)).  For $\hat\tau\in (d/2, d-1)$,
$p_\pm(\hat\tau)\geq 2$
and so we may
let $\tilde {\cal B}_{\tau\pm\varepsilon}$
denote the blow-ups of ${\cal B}_{\tau\pm\varepsilon}$ along
${\Bbb P}_\varepsilon(\Wtau^\pm)$
and $\rho_\pm$ denote
the corresponding blow-down maps.
If $\hat\tau =d-1$, $p_-(d-1)=1$, and ${\Bbb
P}_\varepsilon(\Wtau^-)$ is already a divisor.  In this case we let
$\tilde{\cal B}_{\tau-\varepsilon}={\cal B}_{\tau-\varepsilon}$ and
$\tilde {\cal B}_{\tau+\varepsilon}$ the blow-up of ${\cal
B}_{\tau+\varepsilon}$ along ${\Bbb P}_\varepsilon(\Wtau^+)$, which
may be identified with $\Sigma\times {\cal J}_{d-1}$ by Example 3.5.
In any case, the biholomorphism
$$\sigma_+^{-1}\circ\sigma_-: {\cal B}_{\tau-\varepsilon}\setminus
\sigma_-^{-1}
(\Ztau)\longrightarrow
{\cal B}_{\tau+\varepsilon}\setminus\sigma_+^{-1} (\Ztau)$$
from Corollary 6.2 clearly lifts to a bimeromorphic map
$\tilde\sigma:\tilde{\cal B}_{\tau-\varepsilon}\to\tilde {\cal
B}_{\tau+\varepsilon}$.  On the other hand, the result of [G-S]
proves that $\tilde\sigma$ extends to a continuous bijection on all of
$\tilde {\cal B}_{\tau-\varepsilon}$.  Now the Riemann extension
theorem implies that $\tilde {\cal B}_{\tau-\varepsilon}$ is
biholomorphic to $\tilde {\cal B}_{\tau+\varepsilon}$. We therefore
set $\tilde\Btau= \tilde {\cal B}_{\tau-\varepsilon}
\simeq\tilde {\cal B}_{\tau+\varepsilon}$.

\beginsection 7. Concluding remarks

In [B-D2] we introduced the moduli space of stable pairs of {\it
fixed determinant} (see also [Th]).  These are defined as follows:
Let $\Jd$ denote component of the Jacobian of $\Sigma$ corresponding
to degree $d$ line bundles, and let
$$\det : \Btau\longrightarrow \Jd\leqno(7.1)$$
denote the map $\det(E^{\dbare},\phi)=\Lambda^R E^{\dbare}$.
It was shown in [B-D2] that  det is a holomorphic map of maximal
rank and thus a fibration.  Let $\BtL$ denote the fiber of det over
$L\in \Jd$.  More generally, let
$$\det : \Bhat\longrightarrow \Jd\leqno(7.2)$$
denote the map $\det(E^{\dbare},\phi)=\Lambda^R E^{\dbare}$.
For $L\in \Jd$ let $\Bhat(L)=\det^{-1}(L)$ and
$\Bhat_0(L)=\Bhat(L)\cap \Bhat_0$.
Clearly, $\Bhat(L)$ and  $\Bhat_0(L)$
are preserved by the $\U$ action, and
for any non-critical value of $\hat\tau$,
$f^{-1}(\hat\tau)\bigcap\Bhat(L)/\U$
is biholomorphic to $\BtL$.

It is easily seen that all the constructions performed in the
previous sections commute with the map det and thus one has the
analogous theorems for $\BtL$.  In particular, Theorem 6.1 and
and Corollaries 6.2, 6.4 and 6.5 remain valid by replacing $\Btau$ by
$\BtL$.

Perhaps the most important question is how to resolve the birational
maps of Corollary 6.4.  The problem is that the master space $\Bhat$
is singular along some of the critical sets.  This means that ${\Bbb
P}_\varepsilon(\Wtau^\pm)$, the centers along which we wish to
blow-up, are singular in general.  One way to proceed might be to
desingularize the master space $\Bhat$ as in Kirwan [K2] and extend
the circle action.  However, one would still have to deal with finite
quotient singularities.  Such a description is desirable because by
Corollary 6.5 one would then have a relationship between the moduli
spaces of rank $R$ bundles which is inductive on the rank.  This
could be used to compute, for example, Verlinde dimensions as in
[Th] or perhaps even the cohomology ring structure of these spaces
in a manner similar to [Be-D-W].

\beginsection {} References

\item{[A-B]} Atiyah, M. F. and R. Bott, {\it The Yang-Mills equations
over Riemann surfaces}, Phil. Trans. R. Soc. Lond. A {\bf 308}
 (1982), 523-615.

\item{[B1]} Bradlow, S. B., {\it Special metrics and stability for
holomorphic bundles with global sections}, J. Diff. Geom. {\bf 33}
(1991), 169-214.

\item{[B2]} Bradlow, S. B., {\it Vortices in holomorphic line
bundles over K\"ahler manifolds}, Commun. Math. Phys. {\bf 135}
(1990), 1-18.

\item{[B-D1]} Bradlow, S. B. and G. D. Daskalopoulos, {\it Moduli of
stable pairs for holomorphic bundles over Riemann surfaces}, Int.
J. Math. {\bf 2} (1991), 477-513.

\item{[B-D2]} Bradlow, S. B. and G. D. Daskalopoulos, {\it Moduli of
stable pairs for holomorphic bundles over Riemann surfaces II},
to appear in Int. J. Math.

\item{[Be]} Bertram, A., {\it Stable pairs and stable parabolic
pairs}, preprint, 1992.

\item{[Be-D-W]} Bertram, A., G. Daskalopoulos, and R. Wentworth,
{\it Gromov invariants for holomorphic maps from Riemann surfaces to
Grassmannians}, preprint, 1993.

\item{[D]} Daskalopoulos, G.D., {\it The topology of the space of
stable bundles over a compact Riemann surface}, J. Diff. Geom. {\bf
36} (1992), 699-746.

\item{[G-H]} Griffiths, P. and J. Harris, ``Principles of Algebraic
Geometry", Wiley, New York, 1978.

\item{[G-M]} Goresky, M. and R. MacPherson, {\it On the topology of
torus actions}, in Lecture Notes in Mathematics {\bf 1271},
Springer-Verlag, Berlin, Heidelberg, 1987.

\item{[G-P]} Garcia-Prada, O., {\it Dimensional reduction of stable
bundles, vortices and stable pairs}, preprint.

\item{[G-S]} Guillemin, V. and  S. Sternberg, {\it Birational
equivalence in the symplectic category}, Invent. Math. {\bf 97}
(1989), 485-522.

\item{[K1]} Kirwan, F., ``Cohomology of Quotients in Symplectic and
Algebraic Geometry", \break Princeton University Press, 1984.

\item{[K2]} Kirwan, F., {\it On the homology of compactifications of
moduli spaces of vector bundles over a Riemann surface}, Proc.
London Math. Soc. (3) {\bf 53} (1986), 237-266.

\item{[K-K]} Kaup, L. and B. Kaup, ``Holomorphic Functions of
Several Variables", de Gruyter Studies in Mathematics {\bf 3},
Berlin, New York, 1983.

\item{[Ko]} Kobayashi, S., ``Differential Geometry of Complex Vector
Bundles", Princeton University Press, Princeton, 1987.

\item{[M]} Moishezon, B.G., {\it A criterion for projectivity of
complete algebraic abstract varieties}, A.M.S. Transl. {\bf 63}
(1967), 1-50.

\item{[N]} Newstead, P.E., ``Introduction to Moduli Problems and
Orbit Spaces", Tata Inst. Lectures {\bf 51}, Springer-Verlag,
Heidelberg, 1978.

\item{[S]} Siu, Y.-T., {\it Vanishing theorems for the semipositive
case},  in Lecture Notes in Math. {\bf 1111}, Springer-Verlag,
Heidelberg, 1985, pp. 164-192.

\item{[Th]} Thaddeus, M., {\it Stable pairs, linear systems, and the
Verlinde formula}, preprint, 1992.

\bigskip\bigskip\smalltype
\baselineskip 10pt \normalbaselineskip 10pt
\leftline{Department of Mathematics}
\leftline{University of Illinois}
\leftline{273 Altgeld Hall}
\leftline{1409 W. Green Street}
\leftline{Urbana, IL  61801}
\leftline{email: bradlow@vortex.math.uiuc.edu}
\bigskip
\leftline{Department of Mathematics}
\leftline{Princeton University}
\leftline{Fine Hall}
\leftline{Princeton, NJ 08544}
\leftline{email: daskal@math.princeton.edu}
\bigskip
\leftline{Department of Mathematics}
\leftline{Harvard University}
\leftline{One Oxford Street}
\leftline{Cambridge, MA 02138}
\leftline{email: raw@math.harvard.edu}
\end